\documentclass[a4paper]{amsart}
\usepackage{amssymb} 
\usepackage[T1]{fontenc}
\usepackage[latin1]{inputenc}
\usepackage{amsfonts}
\usepackage{amsxtra}
\usepackage{ae}
\usepackage[all]{xy}
\usepackage{enumerate}

\include{diagram}

\newcommand*{\F}{\mathcal{F}}
\newcommand*{\G}{\mathcal{G}}
\newcommand*{\K}{\mathcal{K}}

\newcommand*{\T}{\mathcal{T}}
\newcommand*{\cotimes}{\hat{\otimes}}

\newcommand*{\DD}{\mathsf{D}}

\newcommand*{\yd}{\mathsf{YD}}
\newcommand*{\AYD}{\mathsf{A}}

\newcommand*{\hit}{\rightharpoonup}
\newcommand*{\hitby}{\leftharpoonup}
\newcommand*{\ihd}{{\raisebox{0cm}[0cm][0cm]{$\scriptstyle{\hat{H}}$}}}
\newcommand*{\red}{\mathsf{r}}

\newcommand*{\cop}{\mathsf{cop}}
\newcommand*{\CH}{\mathbb{C}}

\DeclareMathOperator{\End}{End}

\DeclareMathOperator{\Hom}{Hom}

\DeclareMathOperator{\id}{id}
\DeclareMathOperator{\tr}{tr}

\DeclareMathOperator{\Pol}{Pol}

\newenvironment{bnum}
{\begin{list}{}
    {\setlength{\labelwidth}{15pt}
     \setlength{\leftmargin}{\labelwidth}
    }
}
{\end{list}}

\numberwithin{equation}{section}
\theoremstyle{change}
\newtheorem{theorem}{Theorem}[section]
\newtheorem{prop}[theorem]{Proposition}
\newtheorem{lemma}[theorem]{Lemma}

\newtheorem{definition}[theorem]{Definition}

\begin{document}

\title[Cyclic cohomology and Baaj-Skandalis duality]{Cyclic cohomology and Baaj-Skandalis duality}
\author{Christian Voigt}
\address{School of Mathematics and Statistics \\
         University of Glasgow \\
         University Gardens \\
         Glasgow G12 8QW \\
         United Kingdom 
}
\email{christian.voigt@glasgow.ac.uk}

\subjclass[2000]{19D55}

\maketitle

\begin{abstract}
We construct a duality isomorphism in equivariant periodic cyclic homology analogous to Baaj-Skandalis duality in equivariant Kasparov 
theory. As a consequence we obtain general versions of the Green-Julg theorem and the dual Green-Julg theorem in periodic 
cyclic theory. \\
Throughout we work within the framework of bornological quantum groups, thus in particular incorporating at the same time actions of 
arbitrary classical Lie groups as well as actions of compact or discrete quantum groups. An important ingredient in the construction of our duality 
isomorphism is the notion 
of a modular pair for a bornological quantum group, closely related to the concept introduced by Connes and Moscovici in their work on cyclic 
cohomology for Hopf algebras. 
\end{abstract}

\section{Introduction}

The classical Takesaki-Takai duality theorem for abelian locally compact groups, along with its generalisations to general locally compact groups and 
quantum groups, plays an important role in the study of $ C^* $-dynamical systems. It leads, among many other things, to a duality isomorphism on 
the level of equivariant Kasparov theory due to Baaj and Skandalis \cite{BSKK}, \cite{BSUM}. \\
More precisely, assume that $ S $ and $ \hat{S} $ are the reduced Hopf $ C^* $-algebras associated to a Kac system, 
or more generally, to a regular symmetric multiplicative unitary, see \cite{BSUM}, \cite{Fischerthesis}. Then there exists a canonical isomorphism
$$
J_S: KK^S_*(A, B) \rightarrow KK^{\hat{S}}_*(A \rtimes_\red S, B \rtimes_\red S) 
$$
of the corresponding equivariant Kasparov groups, obtained essentially by forming crossed products of the underlying equivariant Kasparov cycles. 
Here we consider reduced crossed products on the right hand side, equipped with the canonical dual coactions. 
A prototypical example of this situation is that $ S $ and $ \hat{S} $ are associated to a regular locally compact quantum group \cite{KVLCQG}. 
In particular, the above isomorphism holds for locally compact groups and their duals, and it is often in this setting that Baaj-Skandalis duality appears 
in applications, see for instance \cite{Cuntzbc}. \\ 
The main result of this paper is an analogue of the Baaj-Skandalis duality isomorphism in equivariant cyclic homology. Considering actions of, say, locally 
compact quantum groups on $ C^* $-algebras is not suitable in this context. We work with actions of 
bornological quantum groups on bornological algebras instead. This setup unifies various commonly studied 
situations relevant to cyclic homology, in particular it allows to treat actions of discrete groups or quantum groups on 
algebras without further structure in the same way as smooth actions of Lie groups on bornological algebras. 
We show that there exists a natural isomorphism 
$$
J_H: HP^H_*(A, B) \rightarrow HP^{\hat{H}}_*(A \rtimes H, B \rtimes H) 
$$
in equivariant periodic cyclic homology where $ H $ is a bornological quantum group and $ \hat{H} $ is the dual quantum group. Here 
the crossed products on the right hand side are taken in the bornological sense, and equipped with the canonical dual actions. 
This isomorphism is compatible with composition products, and it allows us to derive the Green-Julg theorem in periodic cyclic 
homology for compact quantum groups and its dual for discrete quantum groups. Whereas the Green-Julg theorem in cyclic homology has been studied 
in various setups in the literature \cite{Brylinskibrown}, \cite{Blockthesis}, \cite{Buesthesis}, \cite{AK2}, the dual Green-Julg theorem, 
being considerably more complicated, has not 
received much attention so far. Our approach to the dual Green-Julg theorem in this paper is at the same time more conceptual and more general than the 
computational argument given in \cite{Voigtepch}. \\
The construction of the Baaj-Skandalis duality map in cyclic homology is surprisingly subtle and relies on a detailed analysis of the 
equivariant $ X $-complexes defining the theory. In addition, it requires in an essential way the use of a modular pair for the corresponding bornological 
quantum groups. 
Modular pairs for Hopf algebras, consisting of a group-like element and a character satisfying certain properties, have been originally introduced by 
Connes and Moscovici in their work on cyclic cohomology 
for Hopf algebras \cite{CMtransverse}, \cite{CMcyclichopf}, \cite{CMcyclichopfsymmetry}. Such pairs play an important role in Hopf-cyclic homology 
and cohomology \cite{HKRS}, 
providing coefficients for these theories. In our context, they are used to reformulate the definition of 
equivariant cyclic homology in terms of Yetter-Drinfeld modules instead of anti-Yetter-Drinfeld modules. 
Without this adjustment, the duality map cannot be written down. \\ 
Let us also remark that modular pairs for bornological quantum groups are closely related to the modular elements of the quantum group and its dual. 
More precisely, they provide square roots of these elements. This is completely analogous to the source of the modular character and 
the twisted antipode in the work of Connes and Moscovici \cite{CMtransverse}. At the same time, the notion of a modular pair 
for a bornological quantum group is more restrictive than the purely algebraic concept studied in the setting of general Hopf algebras. \\ 
Let us now describe how the paper is organised. In section \ref{secqg} we review some background material on bornological quantum groups, 
including duality and actions on algebras. Section \ref{secmodularpairs} contains the definition of a modular pair for a bornological 
quantum group, along with some examples. As mentioned above, modular pairs play a crucial role in the identification of 
anti-Yetter-Drinfeld modules and Yetter-Drinfeld modules, which we discuss in section \ref{secydayd}. 
Section \ref{secech} is devoted to a brief review of equivariant cyclic homology. In section \ref{secdualitymap} we define the duality map relating 
equivariant differential forms of an $ H $-algebra
with equivariant differential forms of its crossed product. Finally, in section \ref{secbs} we prove our main result 
and discuss the Green-Julg theorems. \\ 
As already indicated above, we shall work within the monoidal category of complete bornological vector spaces throughout. 
The tensor product $ \cotimes $ in this category is the completed projective bornological tensor product. 
We refer to \cite{Meyerbook} for background information on bornological vector spaces. 
Let us point out that the theory of bornological vector spaces is not needed if one restricts attention to algebraic quantum groups 
in the sense of Van Daele \cite{vDadvances}, acting on algebras without further structure. 
However, the modular properties entering in our discussion are most clearly visible if one includes
examples coming from noncompact Lie groups, and the corresponding bornological quantum groups are outside the scope 
of a purely algebraic framework.

\section{Bornological quantum groups} \label{secqg} 

In this section we review the theory of bornological quantum groups, including duality and their actions on bornological algebras. 
For more information and details we refer to \cite{Voigtbqg}. \\ 
A bornological algebra $ H $ is called essential if the multiplication map induces an isomorphism $ H \cotimes_H H \cong H $. 
The multiplier algebra $ M(H) $ of a bornological algebra $ H $ consists of all two-sided multipliers of $ H $, the latter being 
defined by the usual algebraic conditions. There exists a canonical bounded 
homomorphism $ \iota: H \rightarrow M(H) $. A bounded linear functional $ \phi: H \rightarrow \mathbb{C} $ on a bornological algebra 
is called faithful if $ \phi(st) = 0 $ for all $ t \in H $ implies $ s = 0 $ and $ \phi(st) = 0 $ for all $ s \in H $ implies $ t = 0 $. 
If there exists such a functional the map $ \iota: H \rightarrow M(H) $ is injective, 
and one may view $ H $ as a subset of the multiplier algebra $ M(H) $. \\
In the sequel $ H $ will be an essential bornological algebra with a faithful bounded linear functional. 
For technical reasons we assume moreover that the underlying bornological vector space of $ H $ satisfies the approximation 
property. \\
A left module $ V $ over $ H $ is called essential if the module action $ H \cotimes V \rightarrow V $ induces an 
isomorphism $ H \cotimes_H V \cong V $. A bounded linear map $ f: V \rightarrow W $ between essential $ H $-modules is called $ H $-linear 
if it commutes with the action of $ H $. 
Similarly one defines essential right modules, and we call an algebra homomorphism $ H \rightarrow M(K) $ essential 
if it turns $ K $ into an essential left and right $ H $-module. 
Let $ \Delta: H \rightarrow M(H \cotimes H) $ be an essential homomorphism. The left Galois 
maps $ \gamma_l, \gamma_r: H \cotimes H \rightarrow M(H \cotimes H) $ for $ \Delta $ are defined by 
$$
\gamma_l(s \otimes t) = \Delta(s)(t \otimes 1),\qquad \gamma_r(s \otimes t) = \Delta(s) (1 \otimes t). 
$$
Similarly, the right Galois maps $ \rho_l, \rho_r: H \cotimes H \rightarrow M(H \cotimes H) $ for $ \Delta $ are defined by 
$$
\rho_l(s \otimes t) = (s \otimes 1) \Delta(t), \qquad \rho_r(s \otimes t) = (1 \otimes s) \Delta(t). 
$$
The map $ \Delta $ is called a comultiplication if
$$
(\Delta \cotimes \id) \Delta = (\id \cotimes \Delta)\Delta, 
$$
where both sides are viewed as maps from $ H $ to $ M(H \cotimes H \cotimes H) $. \\ 
Let $ \Delta: H \rightarrow M(H \cotimes H) $ be a comultiplication such that all Galois maps associated 
to $ \Delta $ define bounded linear maps from $ H \cotimes H $ into itself. Then 
a bounded linear functional $ \phi: H \rightarrow \mathbb{C} $ is called left invariant if 
\begin{equation*}
(\id \cotimes \phi)\Delta(t) = \phi(t) 1
\end{equation*}
for all $ t \in H $, where the left hand side has to be interpreted appropriately as a multiplier of $ H $. 
Similarly, one defines right invariant functionals. \\ 
Let us now recall the definition of a bornological quantum group \cite{Voigtbqg}.
\begin{definition} \label{bqgdef}
A bornological quantum group consists of an essential bornological algebra $ H $ satisfying the approximation property, a 
comultiplication $ \Delta: H \rightarrow M(H \cotimes H) $ such that all 
Galois maps associated to $ \Delta $ are isomorphisms, and a faithful left invariant functional $ \phi: H \rightarrow \mathbb{C} $. 
\end{definition}
One can show that the functional $ \phi $ in definition \ref{bqgdef} is unique up to a scalar. 
The definition of a bornological quantum group is equivalent to the definition of an algebraic quantum group in the 
sense of Van Daele \cite{vDadvances} if the underlying bornological vector space carries the fine bornology. \\
In particular, if $ G $ is a compact quantum group then the unital Hopf $ * $-algebra $ \Pol(G) $ of polynomial functions 
on $ G $ can be viewed as a bornological quantum group. 
Further natural examples of bornological quantum groups arise from Lie groups, for instance. More precisely, if $ G $ 
is a possibly noncompact Lie group then the algebra $ C^\infty_c(G) $ of compactly supported smooth functions with the precompact bornology 
is a bornological quantum group. The comultiplication is induced from the group law of $ G $ in this case. \\
The following result from \cite{Voigtbqg} shows that one may view bornological quantum groups as generalised Hopf algebras. 
\begin{theorem} \label{bqchar}
Let $ H $ be a bornological quantum group. Then there exists an essential algebra homomorphism  $ \epsilon: H \rightarrow \mathbb{C} $ and 
a linear isomorphism $ S: H \rightarrow H $ which is both an algebra antihomomorphism and a coalgebra antihomomorphism such that 
\begin{equation*}
(\epsilon \cotimes \id)\Delta = \id = (\id \cotimes \epsilon)\Delta
\end{equation*}
and
\begin{equation*}
\mu(S \cotimes \id) \gamma_r = \epsilon \cotimes \id, \qquad \mu(\id \cotimes S) \rho_l = \id \cotimes \epsilon.
\end{equation*}
In addition, the maps $ \epsilon $ and $ S $ are uniquely determined. 
\end{theorem}
Using the antipode one finds that every bornological quantum groups is equipped with a faithful right invariant functional $ \psi $ 
as well, again unique up to a scalar. We will typically fix the choice of $ \phi $ and $ \psi $ in this way. \\
Due to the existence of invariant functionals one obtains a well-behaved duality theory for bornological quantum groups, extending the 
duality theory of algebraic quantum groups developed by Van Daele. More precisely, let us define bounded linear 
maps $ \F_l, \F_r, \G_l, \G_r $ from $ H $ into the space $ H' $ of bounded linear functionals on $ H $ by 
\begin{align*}
\mathcal{F}_l(t)(r) &= \phi(rt), \qquad \mathcal{F}_r(t)(r) = \phi(tr) \\
\mathcal{G}_l(t)(r) &= \psi(rt), \qquad \mathcal{G}_r(t)(r) = \psi(tr).
\end{align*} 
The images of these maps coincide and determine a vector subspace $ \hat{H} $ of $ H' $. Moreover, there exists a unique bornology 
on $ \hat{H} $ such that these maps are bornological isomorphisms from $ H $ to $ \hat{H} $. Using the transposition of the 
multiplication, comultiplication and counit maps of $ H $, one obtains a canonical bornological quantum group structure on $ \hat{H} $. 
The invariant integrals for $ \hat{H} $ are defined by 
$$
\hat{\phi}(\mathcal{G}_r(t)) = \epsilon(t), \qquad \hat{\psi}(\mathcal{F}_l(t)) = \epsilon(t),
$$
respectively. 
\begin{theorem} 
Let $ H $ be a bornological quantum group. Then $ \hat{H} $ with the structure maps described above is again a 
bornological quantum group. Moreover, the dual of $ \hat{H} $ is canonically isomorphic to $ H $.
\end{theorem}
We will often make use of Sweedler notation in our computations. That is, we write 
$$
\Delta(t) = t_{(1)} \otimes t_{(2)}
$$
for the coproduct of an element $ t $ in a bornological quantum group, and accordingly for higher coproducts. 
This notation has of course only formal meaning. 
Nonetheless, since some of the calculations we have to perform in subsequent sections are quite complicated, 
it will help to organise the arguments efficiently. Let us point out that our computations could be rewritten using 
diagrams of bounded linear maps only, although this would certainly obscure the basic ideas. \\
An $ H $-algebra is by definition an algebra object in the category of essential $ H $-modules. 
We formulate this more explicitly in the following definition.
\begin{definition}
Let $ H $ be a bornological quantum group. An $ H $-algebra is a bornological algebra $ A $ which is at the same 
time an essential $ H $-module such that the multiplication map $ A \cotimes A \rightarrow A $ is $ H $-linear.
\end{definition}
Notice that we do not require the existence of unit elements. We write $ A^+ $ for the unitarisation of 
an $ H $-algebra. As a bornological vector space we have $ A^+ = A \oplus \mathbb{C} $, and multiplication is defined 
in such a way that it becomes a unital algebra, with the copy of $ A $ being contained in $ A^+ $ as an ideal.

\section{Modular pairs} \label{secmodularpairs} 

In this section we discuss the concept of a modular pair for a bornological quantum group. The notion of a modular pair 
was introduced by Connes and Moscovici for arbitrary Hopf algebras in \cite{CMtransverse}, \cite{CMcyclichopf}, \cite{CMcyclichopfsymmetry}, 
see also \cite{HKRS}. We remark that our definition is somewhat more restrictive because it is directly linked to the 
modular elements of the underlying quantum group and its dual. \\
Let $ H $ be a bornological quantum group. Then a group-like element for $ H $ is an invertible element $ \sigma \in M(H) $ such 
that $ \Delta(\sigma) = \sigma \otimes \sigma $ as well as $ S(\sigma) = \sigma^{-1} $ and $ \epsilon(\sigma) = 1 $. 
A character for $ H $ is a nondegenerate algebra 
homomorphism $ \delta $ from $ H $ to $ \mathbb{C} $. The extension of this homomorphism 
to the multiplier algebra will again be denoted by $ \delta $. Using the duality theory of bornological quantum 
groups \cite{Voigtbqg}, it is easy to check that a group-like element for $ H $ is the same thing as a character for the dual 
quantum group $ \hat{H} $, and vice versa. \\
Let us introduce some further notation. 
Given an element $ t \in H $ and $ \omega \in H' $, we define bounded linear functionals $ t \hit \omega $ and $ \omega \hitby t $ on $ H $ by 
$$
(t \hit \omega)(r) = \omega(rt), \qquad (\omega \hitby t)(r) = \omega(tr).
$$
These actions preserve the subspace $ \hat{H} \subset H' $ and turn $ H' $ into a left and right $ H $-module. We will also use 
this in the sequel with the roles of $ H $ and $ \hat{H} $ reversed. 
\begin{definition} \label{modpair}
Let $ H $ be a bornological quantum group and let $ \phi $ be a left invariant Haar functional on $ H $.  
A modular pair $ \tau = (\sigma, \delta) $ for $ H $ consists of a group-like element $ \sigma \in M(H) $ 
and a character $ \delta: H \rightarrow \mathbb{C} $ satisfying the conditions
\begin{bnum}
\item[a)] (Invariance) $ (\phi \cotimes \id)\Delta(t) = \sigma^{-2} \phi(t) $,
\item[b)] (Weak KMS property) $ \phi(rt) = \phi(t \delta \hit(\sigma r \sigma^{-1}) \hitby \delta) $,
\item[c)] (Involutivity) $ S^2(t) = \delta^{-1} \hit (\sigma t \sigma^{-1}) \hitby \delta $, 
\item[d)] (Normalisation) $ \delta(\sigma) = 1 $,
\end{bnum}
where $ r,t \in H $. 
\end{definition}
In the terminology of Connes-Moscovici \cite{CMcyclichopfsymmetry}, our definition corresponds to a modular pair in 
involution. The involutivity condition c) in definition \ref{modpair} plays an important role in the framework of cyclic cohomology for Hopf 
algebras \cite{HKRScyclic}. 
We remark that condition a) in definition \ref{modpair} says that $ \sigma^{-2} $ is the modular element of $ H $. This implies 
in particular that $ \sigma^2 $ is uniquely determined by $ H $, see \cite{Voigtbqg}. We will see below that the same holds true for $ \delta^2 $. \\
Nonetheless, modular pairs for a bornological quantum group are by no means unique. Consider for instance the group 
algebra $ H = \CH[\mathbb{Z}_2] $ of the group $ \mathbb{Z}_2 $. Since in this case both $ H $ and $ \hat{H} $ are 
unimodular we see that $ (1, \epsilon) $ is a modular pair for $ H $. If $ \sigma \in \mathbb{Z}_2 $ denotes the 
generator, viewed as a group-like element of $ H $, then $ (\sigma, \epsilon) $ is easily seen to be a modular pair 
as well. In a similar way we obtain families of nontrivial modular pairs for various other examples of quantum groups. \\
Our first aim is to study how a modular pair transforms under duality. Let $ H $ be a bornological quantum group 
and let $ (\sigma, \delta) $ be a modular pair for $ H $. Then we have 
$$
S^2(t) = \sigma (\delta^{-1} \hit t \hitby \delta) \sigma^{-1}, \qquad S^{-2}(t) = \sigma^{-1} (\delta \hit t \hitby \delta^{-1}) \sigma, 
$$
which implies
$$
S(t) = \sigma (\delta^{-1} \hit S^{-1}(t) \hitby \delta) \sigma^{-1}, \qquad S^{-1}(t) = \sigma^{-1}(\delta \hit S(t) \hitby \delta^{-1}) \sigma
$$
for all $ t \in H $. For an element $ f \in \hat{H} $ we obtain
$$
S^2(f) = \delta (\sigma^{-1} \hit f \hitby \sigma) \delta^{-1}, \qquad S^{-2}(f) = \delta^{-1} (\sigma \hit f \hitby \sigma^{-1}) \delta 
$$
and 
$$
S(f) = \delta (\sigma^{-1} \hit S^{-1}(f) \hitby \sigma) \delta^{-1}, \qquad S^{-1}(f) = \delta^{-1} (\sigma \hit S(f) \hitby \sigma^{-1}) \delta. 
$$
Moreover, if we fix a right invariant integral for $ H $ by $ \psi = \phi S $ then 
\begin{align*}
\psi(rt) &= \phi(S(t) S(r)) = \phi(S(r) \delta \hit (\sigma S(t) \sigma^{-1}) \hitby \delta)\\
&= \phi(S(r) S(\delta^{-1} \hit (\sigma t \sigma^{-1}) \hitby \delta^{-1}))\\
&= \psi(\delta^{-1} \hit (\sigma t \sigma^{-1}) \hitby \delta^{-1} r). 
\end{align*}
This implies 
$$
\psi(tr) = \psi(r \delta \hit (\sigma^{-1} t \sigma) \hitby \delta)
$$
for all $ r,t \in H $. \\
Let us now study the Haar functional $ \hat{\phi} $ on the dual quantum group $ \hat{H} $, normalised as indicated in section 
\ref{secqg}. We compute 
\begin{align*}
(\hat{\phi} \cotimes s) \Delta(\G_r(t)) = \hat{\phi}(\G_r(\delta^{-1} \hit (\sigma s \sigma^{-1}) \hitby \delta^{-1} t)) 
= \delta^{-2}(s) \epsilon(t) 
\end{align*}
for $ s \in H $ and deduce 
$$
(\hat{\phi} \cotimes \id) \Delta(f) = \delta^{-2} \hat{\phi}(f)
$$
for all $ f \in \hat{H} $. Moreover we calculate 
$$
\hat{\phi}(\G_r(s) \G_r(t)) = \psi(t S^{-1}(s)) = \psi(s S(t) \sigma^2)
$$
and 
\begin{align*}
\sigma \hit (\delta \G_r(t) \delta^{-1}) \hitby \sigma &= \delta \G_r(\sigma t \sigma) \delta^{-1} \\
&= \delta^{-1}(t_{(1)}) \G_r(\sigma t_{(2)} \sigma) \delta(t_{(3)})
= \G_r(\delta \hit (\sigma t \sigma) \hitby \delta^{-1})
\end{align*}
using invariance of $ \psi $. Hence we get 
\begin{align*}
\hat{\phi}(\G_r(t) \sigma  \hit (\delta \G_r(s) \delta^{-1}) \hitby \sigma) &= 
\hat{\phi}(\G_r(t) \G_r(\delta \hit (\sigma s \sigma) \hitby \delta^{-1})) \\
&= \psi(t S(\delta \hit (\sigma s \sigma) \hitby \delta^{-1}) \sigma^2) \\
&= \psi(t \delta \hit (\sigma^{-1} S(s) \sigma) \hitby \delta^{-1}) = \psi(t S^{-1}(s))
\end{align*}
and deduce 
$$
\hat{\phi}(fg) = \hat{\phi}(g \sigma \hit (\delta f \delta^{-1}) \hitby \sigma)
$$
for all $ f,g \in \hat{H} $. We have thus proved the following result. 
\begin{prop}\label{modpairdual}
Let $ \tau = (\sigma, \delta) $ be a modular pair for the bornological quantum group $ H $. Then 
$ \hat{\tau} = (\delta, \sigma) $ is a modular pair for the dual quantum group $ \hat{H} $. 
\end{prop}
Let us have look at some examples of modular pairs. If $ H = C^\infty_c(G) $ is the bornological quantum group of compactly supported smooth functions 
on a Lie group $ G $ then we obtain a natural modular pair $ (\Delta^{-\frac{1}{2}}, \epsilon) $ for $ H $ where 
$ \Delta \in C^\infty(G) $ is the modular function of $ G $. \\
The algebra of polynomial functions $ \Pol(G) $ on a compact quantum group $ G $, viewed as a bornological quantum
group with the fine bornology, is equipped with a canonical modular 
pair $ (1, f_{1/2}) $, where $ (f_z)_{z \in \mathbb{C}} $ denotes the family of Woronowicz characters, see \cite{KS}. 
Dually, using proposition \ref{modpairdual} it follows immediately that the algebra of finitely supported functions on a discrete quantum group 
is equipped with a canonical modular pair as well. \\
Finally, modular pairs are compatible with various natural constructions. For instance, 
if $ H_1 $ and $ H_2 $ are bornological quantum groups equipped with modular pairs $ (\sigma_1, \delta_1) $ and $ (\sigma_2, \delta_2) $, 
respectively, then $ (\sigma_1 \cotimes \sigma_2, \delta_1 \cotimes \delta_2)  $ is a modular 
pair for the tensor product $ H_1 \cotimes H_2 $.

\section{Anti-Yetter-Drinfeld modules and Yetter-Drinfeld modules} \label{secydayd}

In this section we discuss anti-Yetter-Drinfeld modules and Yetter-Drinfeld modules for a bornological quantum group and relate them using 
modular pairs. In the algebraic context of Hopf algebras these results are discussed, in a more general setup, in \cite{PSgeneralizedayd}. \\
We begin with the definition of an anti-Yetter-Drinfeld module, compare \cite{HKRS}, \cite{Voigtechq}. 
\begin{definition} \label{defayd}
Let $ H $ be a bornological quantum group. An $ H $-anti-Yetter-Drinfeld module is a smooth left $ H $-module $ M $ 
which is also a smooth left $ \hat{H} $-module such that 
$$
t \cdot (f \cdot m) = (S^2(t_{(1)}) \hit f \hitby S^{-1}(t_{(3)})) \cdot (t_{(2)} \cdot m) 
$$
for all $ t \in H, f \in \hat{H} $ and $ m \in M $.
\end{definition}
A bounded linear map $ f: M \rightarrow N $ between $ H $-anti-Yetter-Drinfeld modules is called a morphism of anti-Yetter-Drinfeld modules 
if it is $ H $-linear and $ \hat{H} $-linear. We write $ \Hom_{\AYD(H)}(M,N) $ for the space of all such morphisms. \\
Anti-Yetter-Drinfeld modules are a basic ingredient in the construction of equivariant cyclic homology. 
We remark that every anti-Yetter-Drinfeld module $ M $ is equipped with a canonical symmetry operator $ T: M \rightarrow M $, 
see \cite{Voigtechq}.  
Definition \ref{defayd} is a variant of the following, more widely known concept of a Yetter-Drinfeld module. 
\begin{definition} 
Let $ H $ be a bornological quantum group. An $ H $-Yetter-Drinfeld module is a smooth left $ H $-module $ M $ which is also a smooth 
left $ \hat{H} $-module such that 
\begin{equation*}
t \cdot (f \cdot m) = (t_{(1)} \hit f \hitby S^{-1}(t_{(3)})) \cdot (t_{(2)} \cdot m) 
\end{equation*}
for all $ t \in H, f \in \hat{H} $ and $ m \in M  $. 
\end{definition}
Thus, if $ S^2 = \id $ we see that anti-Yetter-Drinfeld modules are the same thing as Yetter-Drinfeld modules. 
A bounded linear map $ f: M \rightarrow N $ between $ H $-Yetter-Drinfeld modules is called a morphism of Yetter-Drinfeld modules 
if it is $ H $-linear and $ \hat{H} $-linear. We write $ \Hom_{\DD(H)}(M,N) $ for the space of all such morphisms. \\
It is easy to check that the Yetter-Drinfeld compatibility condition is self-dual. That is, a smooth 
left $ H $-module which is also a smooth left $ \hat{H} $-module is an $ H $-Yetter-Drinfeld 
module iff it is an $ \hat{H} $-Yetter-Drinfeld module. The corresponding statement for anti-Yetter-Drinfeld modules does 
not hold in general. \\
Now assume that $ (\sigma, \delta) $ is a modular pair for the bornological quantum group $ H $. 
Then $ \mathbb{C} $ becomes an $ H $-anti-Yetter-Drinfeld module using the actions determined on $ 1 \in \mathbb{C} $ by
$$
t \bullet 1 = \delta(t), \qquad f \bullet 1 = \sigma^{-1}(f)
$$
for $ t \in H $ and $ \sigma \in \hat{H} $, compare \cite{HKRScyclic}. Indeed, we have 
\begin{align*} 
t \bullet f \bullet 1 &= \sigma^{-1}(f) \delta(t) \\
&= f(S^{-1}(t_{(2)}) t_{(1)} \hitby \delta \sigma^{-1}) \\
&= f(S^{-1}(t_{(2)}) \sigma^{-1} \delta \hit S^2(t_{(1)})) \\
&= f_{(1)}(S^{-1}(t_{(3)})) f_{(3)}(S^2(t_{(1)})) \sigma^{-1}(f_{(2)}) \delta(t_{(2)}) \\
&= (S^2(t_{(1)}) \hit f \hitby S^{-1}(t_{(3)})) \bullet (t_{(2)} \bullet 1) 
\end{align*} 
for $ f \in H $ and $ t \in \hat{H} $. \\ 
This anti-Yetter-Drinfeld module structure on $ \mathbb{C} $ can be used to transform $ H $-anti-Yetter-Drinfeld modules into $ H $-Yetter-Drinfeld modules 
and vice versa. More precisely, assume that $ M $ is an $ H $-anti-Yetter-Drinfeld module, and let us denote 
the action of $ f \in H $ and $ t \in \hat{H} $ on $ m \in M $ by $ f \cdot m $ and $ t \cdot m $,
respectively. We define new actions of $ H $ and $ \hat{H} $ on $ M $ by 
$$
t \bullet m = \delta^{-1}(t_{(1)}) t_{(2)} \cdot m, \qquad f \bullet m = \sigma(f_{(2)}) f_{(1)} \cdot m
$$
for $ t \in H $ and $ f \in \hat{H} $, respectively.
Then we compute 
\begin{align*}
t \bullet (f \bullet m) &= \delta^{-1}(t_{(1)}) \sigma(f_{(2)}) t_{(2)} \cdot (f_{(1)} \cdot m) \\
&= \delta^{-1}(t_{(1)}) \sigma(f_{(2)}) (S^2(t_{(2)}) \hit f_{(1)} \hitby S^{-1}(t_{(4)})) \cdot (t_{(3)} \cdot m) \\
&= f_{(1)}(S^{-1}(t_{(3)})) f_{(3)}(S^2(t_{(1)}) \hitby \delta^{-1} \sigma) f_{(2)} \cdot (t_{(2)} \cdot m) \\
&= f_{(1)}(S^{-1}(t_{(3)})) f_{(3)}(\sigma \delta^{-1} \hit t_{(1)}) f_{(2)} \cdot (t_{(2)} \cdot m) \\
&= (t_{(1)} \hit f \hitby S^{-1}(t_{(3)})) \bullet (t_{(2)} \bullet m) 
\end{align*}
and deduce that $ M $ becomes a Yetter-Drinfeld module in this way. 
\begin{lemma}
Let $ H $ be a bornological quantum group equipped with a modular pair. 
The above construction defines an isomorphism of categories between the category of $ H $-anti-Yetter-Drinfeld modules and 
the category of $ H $-Yetter-Drinfeld modules. 
\end{lemma}
\proof Assume that $ \phi: M \rightarrow N $ is a morphism of $ H $-anti-Yetter-Drinfeld modules, that is, 
an $ H $-linear and $ \hat{H} $-linear bounded linear map. Then for $ t \in H $ we compute 
$$
\phi(t \bullet m) = \phi(\delta^{-1}(t_{(1)}) t_{(2)} \cdot m) = \delta^{-1}(t_{(1)}) t_{(2)} \cdot \phi(m) = t \bullet \phi(m),
$$
and obtain similarly $ \phi(f \bullet m) = f \bullet \phi(m) $ for $ f \in \hat{H} $. 
It follows that $ \phi $ is also a morphism of $ H $-Yetter-Drinfeld modules
with respect to the structure maps defined above. The remaining assertions are obvious. \qed

\section{Equivariant periodic cyclic homology} \label{secech}

In this section we review the definition of equivariant periodic cyclic homology for bornological quantum groups \cite{Voigtechq}. 
Moreover, the canonical constructions in terms of anti-Yetter-Drinfeld module picture will be translated into the framework of Yetter-Drinfeld modules, 
using modular pairs and results from the previous section. \\
Let $ H $ be a bornological quantum group. If $ A $ is an $ H $-algebra we obtain a left action of $ H $ on the space $ H \cotimes \Omega^n(A) $ by defining
\begin{equation*}
r \cdot (t \otimes \omega) = r_{(3)} t S(r_{(1)}) \otimes r_{(2)}
\cdot \omega
\end{equation*}
for $ r,t \in H $ and $ \omega \in \Omega^n(A) $. Here $ \Omega^0(A) = A $ and $ \Omega^n(A) = A^+ \cotimes A^{\cotimes n} $ for $ n > 0 $, 
and we recall that $ A^+ $ denotes the unitarisation of $ A $. 
Moreover there is a left action of the dual quantum group $ \hat{H} $ given by 
\begin{equation*}
f \cdot (t \otimes \omega) = f(t_{(2)}) t_{(1)} \otimes \omega.
\end{equation*}
The equivariant $ n $-forms $ \Omega^n_H(A) $ are defined to be the space $ H \cotimes \Omega^n(A) $ together with the $ H $-action and the $ H $-coaction 
described above. It is straightforward to check that $ \Omega^n_H(A) $ becomes an $ H $-anti-Yetter-Drinfeld module in this way. 
We write $ \Omega_H(A) $ for the direct sum of the spaces $ \Omega^n_H(A) $. \\
Let us define operators $ d $ and $ b_H $ on $ \Omega_H(A) $ by
\begin{equation*}
d(t \otimes \omega) = t \otimes d\omega
\end{equation*}
and
\begin{equation*}
b_H(t \otimes \omega da) = (-1)^{|\omega|} (t \otimes \omega a -
t_{(2)} \otimes (S^{-1}(t_{(1)}) \cdot a) \omega).
\end{equation*}
The operator $ b_H $ should is a twisted version of the usual Hochschild boundary, and it satisfies $ b_H^2 = 0 $ as in the nonequivariant situation. 
Explicitly, we have
\begin{align*}
b_H(t \otimes a^0 da^1 \cdots da^n) &= \sum_{j = 0}^n d^j_H(t \otimes a^0 da^1 \cdots da^n)
\end{align*}
where 
\begin{align*}
d^0_H(t \otimes a^0 da^1 \cdots da^n) &= t \otimes a^0 a^1 da^2 \cdots da^n \\
d^j_H(t \otimes a^0 da^1 \cdots da^n) &= (-1)^j t \otimes a^0 da^1 \cdots d(a^j a^{j + 1}) \cdots da^n \qquad \quad \text{for}\; 0 < j < n \\
d^n_H(t \otimes a^0 da^1 \cdots da^n) &= (-1)^n t_{(2)} \otimes (S^{-1}(t_{(1)}) \cdot a^n) a^0 da^1 \cdots da^{n - 1}.
\end{align*}
In addition, we define the equivariant Connes operator $ B_H $ on $ \Omega^n_H(A) $ by 
\begin{equation*}
B_H(t \otimes a^0da^1 \cdots da^n) = \sum_{i = 0}^n (-1)^{ni}
t_{(2)} \otimes S^{-1}(t_{(1)}) \cdot(da^{n + 1 - i} \cdots da^n)
da^0 \cdots da^{n - i}. 
\end{equation*}
Finally, the canonical symmetry operator $ T $ on $ \Omega_H(A) $ arising from the anti-Yetter-Drinfeld module structure is given by
\begin{equation*}
T(t \otimes \omega) = t_{(2)} \otimes S^{-1}(t_{(1)}) \cdot \omega
= S^{-1}(t_{(2)}) \cdot (t_{(1)} \otimes \omega).
\end{equation*}
All operators constructed so far are morphisms of anti-Yetter-Drinfeld modules. 
\begin{prop} 
Let $ H $ be a bornological quantum group and let $ A $ be an $ H $-algebra. The space $ \Omega_H(A) $ of equivariant differential 
forms is a paramixed complex in the category of anti-Yetter-Drinfeld modules, that is, the relations 
$$
b_H^2 = 0 = B_H^2, \qquad B_H b_H + b_H B_H = \id - T 
$$
hold on $ \Omega_H(A) $. 
\end{prop}
For the definition of equivariant cyclic homology in the Cuntz-Quillen picture only a small part of this paramixed complex 
is needed. 
\begin{definition} Let $ A $ be an $ H $-algebra. The equivariant $ X $-complex 
$ X_H(A) $ of $ A $ is the paracomplex 
\begin{equation*}
    \xymatrix{
      {X_H(A) \colon \ }
        {\Omega^0_H(A)\;} \ar@<1ex>@{->}[r]^-{d} &
          {\;\Omega^1_H(A)/ b_H(\Omega^2_H(A)).} 
            \ar@<1ex>@{->}[l]^-{b_H} 
               } 
\end{equation*}
\end{definition}
We are interested in the equivariant $ X $-complex of the periodic tensor algebra $ \mathcal{T}A $ of an 
$ H $-algebra $ A $. Roughly speaking, the periodic tensor algebra $ \mathcal{T}A $ is the pro-$ H $-algebra 
obtained from the usual tensor algebra of $ A $ by a formal completion procedure. All the above 
constructions for $ H $-algebras carry over to pro-$ H $-algebras, we refer to \cite{Voigtechq} for 
the details. \\
We also recall from \cite{Voigtechq} that for every $ H $-algebra $ A $ one may form the associated crossed product $ A \rtimes H $. The 
underlying bornological vector space of $ A \rtimes H $ is $ A \cotimes H $, and we write $ a \rtimes r $ for a simple tensor 
in $ A \rtimes H $. Using Sweedler notation, the multiplication is defined by 
\begin{equation*}
(a \rtimes r)(b \rtimes t) = a r_{(1)} \cdot b \rtimes r_{(2)}t
\end{equation*}
for $ a, b \in A $ and $ r,t \in H $. Moreover, on $ A \rtimes H $ one has the dual action of $ \hat{H} $ defined by 
$$
f \cdot (a \rtimes r) = a \rtimes f \hit r,  
$$
and the corresponding double crossed product $ A \rtimes H \rtimes \hat{H} $ is naturally an $ H $-algebra again. \\ 
Let us now define the equivariant periodic cyclic homology groups. 
\begin{definition} Let $ H $ be a bornological quantum group and let 
$ A $ and $ B $ be $ H $-algebras. 
The equivariant periodic cyclic homology of $ A $ and $ B $ is 
\begin{equation*}
HP^H_*(A,B) = 
H_*(\Hom_{\AYD(H)}(X_H(\mathcal{T}(A \rtimes H \rtimes \hat{H})), 
X_H(\mathcal{T}(B \rtimes H \rtimes \hat{H})))
\end{equation*}
\end{definition}
We consider the usual differential for a $ \Hom $-complex on the right hand side of this definition, and although 
the equivariant $ X $-complexes are only paracomplexes, this $ \Hom $-complex is indeed a complex in the usual sense. 
As in the group case $ HP^H_* $ is a bifunctor, contravariant in the first variable and
covariant in the second variable. We define 
$$ 
HP^H_*(A) = HP^H_*(\mathbb{C}, A), \qquad HP_H^*(A) = HP^H_*(A, \mathbb{C}) 
$$ 
and call these groups the equivariant periodic cyclic homology and cohomology of $ A $, respectively. There is a natural associative product
\begin{equation*}
HP^H_i(A,B) \times HP^H_j(B,C) \rightarrow HP^H_{i + j}(A,C),
\qquad (x,y) \mapsto y \circ x
\end{equation*}
induced by the composition of maps. 
For further general properties of $ HP^H_* $ we refer to \cite{Voigtechq}. \\
Let us now assume that $ (\sigma, \delta) $ is a modular pair for the bornological quantum group $ H $. We shall describe explicitly 
the Yetter-Drinfeld module structure on $ \Omega_H(A) $ obtained using the constructions from section \ref{secydayd}. 
Firstly, we observe that the structure maps of $ \Omega_H(A) = H \cotimes \Omega(A) $, viewed as a Yetter-Drinfeld module, are given by 
$$
r \bullet (t \otimes \omega) = \delta^{-1}(r_{(1)}) r_{(4)} t S(r_{(2)}) \otimes r_{(3)} \cdot \omega
$$
and 
$$ 
f \bullet (t \otimes \omega) = \sigma(f_{(2)}) f_{(1)}(t_{(2)}) t_{(1)} \otimes \omega
$$
for $ r \in H $ and $ f \in \hat{H} $, respectively. \\ 
It will be convenient to work with a slightly different description of $ \Omega_H(A) $. More precisely, consider 
$ \Omega^{\yd}_H(A) = \hat{H} \cotimes \Omega(A) $ equipped with the actions 
$$
r \bullet (g \otimes \omega) = r_{(1)} \hit g \hitby S^{-1}(r_{(4)}) \otimes r_{(2)} \delta(r_{(3)}) \cdot \omega.
$$
for $ r \in H $ and 
$$ 
f \bullet (g \otimes \omega) = fg \otimes \omega
$$ 
for $ f \in \hat{H} $. We define boundary operators $ b_H $ and $ B_H $ on $ \Omega^{\yd}_H(A) $ by 
$$
b_H(g \otimes \omega da) = (-1)^{|\omega|} (g \otimes \omega a - gS^{-1}(a_{(1)}) \hitby \sigma \otimes a_{(0)} \omega 
$$
and 
\begin{align*}
B_H(g \otimes &a^0da^1 \cdots da^n) \\
&= \sum_{i = 0}^n (-1)^{ni}
g S^{-1}(a^{n + 1 - i}_{(1)} \cdots a^n_{(1)}) \hitby \sigma \otimes da^{n + 1 - i}_{(0)} \cdots da^n_{(0)} 
da^0 \cdots da^{n - i}, 
\end{align*}
respectively. In these formulas we use Sweedler notation for the coaction of $ \hat{H} $ on $ A $ corresponding to the given $ H $-algebra structure, 
compare \cite{Voigtbqg}. \\ 
It is not hard to check directly that $ \Omega_H^{\yd} $ is a paramixed complex of Yetter-Drinfeld modules. Alternatively, 
this is a consequence of the following result. 
\begin{lemma} \label{reformlemma}
Let $ H $ be a bornological quantum group equipped with a modular pair $ (\sigma, \delta) $. If $ A $ is an $ H $-algebra then the bounded 
linear map $ \lambda: \Omega^{\yd}_H(A) \rightarrow \Omega_H(A) $ defined by 
$$
\lambda(f \otimes \omega)(h) = \hat{\psi}(h \sigma \hit f \delta) \otimes \omega
$$
is an isomorphism of paramixed complexes of $ H $-Yetter-Drinfeld modules. 
\end{lemma} 
\proof It is clear that $ \lambda $ is an isomorphism of bornological vector spaces. We calculate 
\begin{align*}
r \bullet &\lambda(f \otimes \omega)(g) = \delta^{-1}(r_{(1)}) r_{(4)}(g_{(1)}) \hat{\psi}(g_{(2)} \sigma \hit f \delta) S(r_{(2)})(g_{(3)}) 
\otimes r_{(3)} \cdot \omega \\
&= \delta^{-1}(r_{(1)}) r_{(4)}(g_{(1)}) \hat{\psi}(g_{(2)} f_{(1)} \delta) S(r_{(2)})(g_{(3)}f_{(2)} \delta \delta^{-1} S(f_{(3)})) \sigma(f_{(4)}) 
\otimes r_{(3)} \cdot \omega \\
&= \delta^{-1}(r_{(1)}) r_{(4)}(g_{(1)}) \hat{\psi}(g_{(2)} f_{(1)} \delta) S(r_{(2)})(\delta^{-1} S(f_{(2)})) \sigma(f_{(3)}) 
\otimes r_{(3)} \cdot \omega \\
&= \delta^{-1}(r_{(1)}) r_{(4)}(g_{(1)}) \hat{\psi}(g_{(2)} f_{(1)} \delta) S(r_{(2)})(\delta^{-1} S(f_{(2)}) \hitby \sigma^{-1})  
\otimes r_{(3)} \cdot \omega \\
&= \delta^{-1}(r_{(1)}) r_{(4)}(g_{(1)}) \hat{\psi}(g_{(2)} f_{(1)} \delta) S(r_{(2)})(\sigma^{-1} \hit S^{-1}(f_{(2)}) \delta^{-1})
\otimes r_{(3)} \cdot \omega \\
&= \delta^{-1}(r_{(1)}) r_{(4)}(g_{(1)}) \hat{\psi}(g_{(2)} f_{(1)} \delta) r_{(2)}(\delta f_{(2)} \hitby \sigma)
\otimes r_{(3)} \cdot \omega \\
&= r_{(3)}(g_{(1)} f_{(2)} \delta \delta^{-1} S^{-1}(f_{(1)})) \hat{\psi}(g_{(2)} f_{(3)} \delta) r_{(1)}(f_{(4)} \hitby \sigma)
\otimes r_{(2)} \cdot \omega \\
&= r_{(3)}(\delta^2 \delta^{-1} S^{-1}(f_{(1)})) \hat{\psi}(g \sigma \hit f_{(2)} \delta) r_{(1)}(f_{(3)}) \otimes r_{(2)} \cdot \omega \\
&= r_{(4)}(S^{-1}(f_{(1)})) \hat{\psi}(g \sigma \hit f_{(2)} \delta) r_{(1)}(f_{(3)}) \otimes r_{(2)} \delta(r_{(3)}) \cdot \omega \\
&= \lambda(r \bullet (f \otimes \omega))(g)
\end{align*}
as desired. Moreover we have 
$$
(g \bullet \lambda(f \otimes \omega))(h) = \hat{\psi}(h \sigma \hit g \sigma \hit f \delta) \otimes \omega
= \hat{\psi}(h (\sigma \hit gf) \delta) = \lambda(gf \otimes \omega)(h), 
$$
which means that $ \lambda $ is $ \hat{H} $-linear. \\ 
Let us now consider the boundary operator $ b_H $. We have 
\begin{align*}
b_H &\lambda(f \otimes \omega da) = b_H(\G_l(\sigma \hit f \delta) \otimes \omega da) \\
&= (-1)^{|\omega|} (\G_l(\sigma \hit f \delta) \otimes \omega a) - 
\G_l(\sigma \hit f \delta \delta^{-1} \sigma \hit S^{-1}(a_{(1)}) \hitby \sigma \delta) \otimes a_{(0)} \omega) \\
&= (-1)^{|\omega|} (\G_l(\sigma \hit f \delta) \otimes \omega a) - 
\G_l(\sigma \hit (f S^{-1}(a_{(1)}) \hitby \sigma) \delta) \otimes a_{(0)} \omega) \\
&= \lambda b_H(f \otimes \omega da). 
\end{align*}
Since $ \lambda $ clearly commutes with the exterior differential operator $ d $ one checks that $ \lambda $ 
is compatible with the Connes operators as well, in the sense that $ \lambda B_H = B_H \lambda $. This finishes 
the proof. \qed \\ 
In the sequel we will consider $ \Omega_H(A) $ as a Yetter-Drinfeld module with the above actions and drop 
reference to the modular pair in our notation. Using lemma \ref{reformlemma} we will also switch tacitly from $ \Omega_H(A) $ to $ \Omega^{\yd}_H(A) $ whenever 
this is convenient. \\
We recall from \cite{Voigtechq} that the double crossed product $ A \rtimes H \rtimes \hat{H} $ of an $ H $-algebra $ A $ 
is $ H $-equivariantly isomorphic to $ A \cotimes \K_H $, where $ \K_H $ is the algebra of smooth kernels on $ H $. 
This is the bornological version of Takesaki-Takai duality. 
Here we note that we may identify $ \K_H = H \cotimes \hat{H} $ equipped with the multiplication 
$$
(r \otimes f)(s \otimes g) = f(s) r \otimes g.  
$$
Observe moreover that $ \K_H $ is naturally a subalgebra of the algebra $ \End(H) $ of bounded linear endomorphisms of $ H $. 
In the presence of modular pairs, it will be convenient to work with a very specific map implementing the Takesaki-Takai isomorphism. 
\begin{prop} \label{TTduality}
Let $ H $ be a bornological quantum group equipped with a modular pair $ (\sigma, \delta) $ and let $ A $ be an $ H $-algebra. 
Then the map $ \gamma_A: A \rtimes H \rtimes \hat{H} \rightarrow A \cotimes \mathcal{K}_H $ 
$$
\gamma_A(a \rtimes r \rtimes f) = (\hat{\psi}_{(3)} \delta \hit S(r_{(1)}) \sigma) \cdot a \otimes \hat{\psi}_{(2)} S(r_{(2)}) \otimes 
S^{-1}(f)(\hat{\psi}_{(1)}) \G_l(S^{-1}(\hat{\psi}_{(4)}))
$$
is an $ H $-equivariant algebra isomorphism. 
\end{prop}
\proof This is completely analogous to proposition 3.7 in \cite{Voigtechq}, but we include the argument for the sake of completeness. 
Observe first that the right invariant functional on $ \hat{H} $ in this formula is diagonalised formally by using the transposition of multiplication 
in $ \hat{H} $. 
Evaluating on an element $ x \in H $ yields
\begin{equation*}
\gamma_A(a \rtimes r \rtimes f)(x) = (x_{(3)} \delta \hit S(r_{(1)}) \sigma) \cdot a \otimes x_{(2)} S(r_{(2)}) \otimes S^{-1}(f)(x_{(1)}). 
\end{equation*}
Using this description of $ \gamma_A $ we compute 
\begin{align*}
&\gamma_A((a^1 \rtimes r^1 \rtimes f^1)(a^2 \rtimes r^2 \rtimes f^2))(x) = 
\gamma_A((a^1 \rtimes r^1)(a^2 \rtimes r^2_{(1)}) \rtimes f^1_{(1)}(r^2_{(2)}) f^1_{(2)} f^2)(x) \\
&= \gamma_A(a^1 r^1_{(1)} \cdot a^2 \rtimes r^1_{(2)} r^2_{(1)} \rtimes f^1_{(1)}(r^2_{(2)}) f^1_{(2)} f^2)(x)  \\
&= x_{(3)} \delta \hit S(r^1_{(2)} r^2_{(1)}) \sigma \cdot (a^1 r^1_{(1)} \cdot a^2) \otimes x_{(2)} S(r^1_{(3)} r^2_{(2)}) \otimes 
S^{-1}(f^1_{(2)} f^2)(x_{(1)}) f^1_{(1)}(r^2_{(3)}) \\
&= (x_{(3)} \delta \hit S(r^2_{(1)}) \sigma S^{-1}(r^1_{(2)}) \hitby \delta) \cdot (a^1 r^1_{(1)} \cdot a^2) \otimes 
x_{(2)} S(r^1_{(3)} r^2_{(2)}) \otimes 
S^{-1}(f^1_{(2)} f^2)(x_{(1)}) f^1_{(1)}(r^2_{(3)}) \\
&= x_{(3)} \delta \hit S(r^2_{(1)}) (\delta \hit S(r^1_{(1)}) \sigma \cdot a^1 \sigma \cdot a^2) \otimes \\
&\qquad x_{(2)} S(r^1_{(2)} r^2_{(2)}) \otimes 
S^{-1}(f^1_{(2)} f^2)(x_{(1)}) f^1_{(1)}(r^2_{(3)}) \\
&= (x_{(4)}S(r^2_{(2)}) \delta \hit S(r^1_{(1)}) \sigma \cdot a^1)(x_{(5)} \delta \hit S(r^2_{(1)}) \sigma \cdot a^2) \otimes \\
&\qquad x_{(3)} S(r^2_{(3)}) S(r^1_{(2)}) \otimes S^{-1}(f^2)(x_{(1)}) S^{-1}(f^1)(x_{(2)} S(r^2_{(4)})) \\
&= \gamma_A(a^1 \rtimes r^1 \rtimes f^1)(x_{(2)}S(r^2_{(2)})) S^{-1}(f^2)(x_{(1)}) x_{(3)} \delta \hit S(r^2_{(1)}) \sigma \cdot a^2  \\
&= \gamma_A(a^1 \rtimes r^1 \rtimes f^1) \gamma_A(a^2 \rtimes r^2 \rtimes f^2)(x). 
\end{align*}
Similarly, one checks 
\begin{align*}
&\gamma_A(t \cdot (a \rtimes r \rtimes f))(x) = \gamma_A(a \rtimes r \rtimes f_{(1)})(x) f_{(2)}(t) \\
&= (x_{(3)} \delta \hit S(r_{(1)}) \sigma) \cdot a \otimes x_{(2)} S(r_{(2)}) \otimes S^{-1}(f_{(1)})(x_{(1)}) f_{(2)}(t) \\
&= (x_{(3)} \delta \hit S(r_{(1)}) \sigma) \cdot a \otimes x_{(2)} S(r_{(2)}) \otimes S^{-1}(f)(S(t) x_{(1)}) \\
&= (t_{(2)} S(t_{(3)}) x_{(3)} \delta \hit S(r_{(1)}) \sigma) \cdot a \otimes t_{(1)} S(t_{(4)}) x_{(2)} S(r_{(2)}) 
\otimes S^{-1}(f)(S(t_{(5)}) x_{(1)}) \\
&= t \cdot \gamma_A(a \rtimes r \rtimes f)(x)
\end{align*}
where we consider the natural $ H $-action 
\begin{equation*}
t \cdot (a \otimes T)(x) = t_{(2)} \cdot a \otimes t_{(1)} T(S(t_{(3)})x).
\end{equation*}
on $ A \otimes \End(H) $, and the corresponding induced action on $ A \otimes \mathcal{K}_H $. 
Finally, it is not hard to verify that $ \gamma_A $ is a bornological isomorphism. \qed

\section{The duality map} \label{secdualitymap}

In this section we define a chain map relating $ H $-equivariant differential forms of an $ H $-algebra $ A $ 
with $ \hat{H} $-equivariant differential forms of the crossed product $ A \rtimes H $. This duality map is then 
analysed carefully. \\ 
Throughout we assume that $ H $ is a bornological quantum group equipped with a fixed modular pair $ (\sigma, \delta) $. We shall use the 
twisted Fourier transform $ \F: \hat{H} \rightarrow H $ given by 
$$
\hat{\F}(f)(h) = \hat{\psi}(S^{-1}(f) \sigma \hit h \delta).
$$
It is easy to check that $ \hat{\F} $ is a bornological isomorphism. \\ 
Let $ A $ be an $ \hat{H} $-algebra. We define a bounded linear map 
$ \tau_H: \Omega_H(A \rtimes \hat{H}) \rightarrow \Omega_\ihd(A) $ by the formula
\begin{align*}
\tau_H&(f \otimes (a^0 \rtimes g^0)d(a^1 \rtimes g^1) \cdots d(a^n \rtimes g^n)) \\
&= \hat{\F}(f_{(n + 3)} g^0_{(2)} g^1_{(3)} \cdots g^n_{(n + 2)} \sigma \hit S(f_{(1)})) \otimes 
(f_{(2)}S^{-1}(g^0_{(1)} g^1_{(2)} \cdots g^n_{(n + 1)}) \cdot  a^0) \\
&\qquad d(f_{(3)} S^{-1}(g^1_{(1)} \cdots g^n_{(n)})\cdot a^1) \cdots 
d(f_{(n + 2)} S^{-1}(g^n_{(1)}) \cdot a^n). 
\end{align*}
Remark that both spaces $ \Omega_H(A \rtimes \hat{H}) $ and $ \Omega_\ihd(A) $ are $ H $-Yetter-Drinfeld modules. 
\begin{prop} \label{duality1}
The map $ \tau_H: \Omega_H(A \rtimes \hat{H}) \rightarrow \Omega_\ihd(A) $ is a map of Yetter-Drinfeld modules and a chain map with respect 
to the equivariant Hochschild and Connes boundary operators. 
\end{prop}
\proof For $ t \in H $ we calculate
\begin{align*}
&\tau_H(t\bullet (f \otimes (a^0 \rtimes g^0) d(a^1 \rtimes g^1) \cdots d(a^n \rtimes g^n)))(k) \\
&= \tau_H(f_{(1)}(S^{-1}(t_{(4)})) f_{(2)} f_{(3)}(t_{(1)}) \otimes \\
&\qquad t_{(2)} \delta(t_{(3)}) \cdot ((a^0 \rtimes g^0) d(a^1 \rtimes g^1) \cdots d(a^n \rtimes g^n)))(k) \\
&= \tau_H(f_{(1)}(S^{-1}(t_{(4)})) f_{(2)} f_{(3)} (t_{(1)}) \otimes (g^0_{(2)} \cdots g^n_{(2)})(t_{(2)}) \delta(t_{(3)})\\
&\qquad (a^0 \rtimes g^0_{(1)}) d(a^1 \rtimes g^1_{(1)}) \cdots d(a^n \rtimes g^n_{(1)})))(k) \\
&= \hat{\F}(t_{(1)} \hit f_{(n + 4)} t_{(2)} \hit (g^0_{(2)} g^1_{(3)} \cdots g^n_{(n + 2)}) 
t_{(3)}(\delta S^{-1}(f_{(1)})) \sigma \hit S(f_{(2)}))(k) \otimes \\
&\qquad (f_{(3)}S^{-1}(g^0_{(1)} \cdots g^n_{(n + 1)}) \cdot  a^0) d(f_{(4)} S^{-1}(g^1_{(1)} \cdots g^n_{(n)})\cdot a^1) \cdots \\
&\qquad \cdots d(f_{(n + 3)} S^{-1}(g^n_{(1)}) \cdot a^n), 
\end{align*}
and by the modular properties of the antipode this is equal to
\begin{align*}
&\quad \hat{\F}(t_{(1)} \hit f_{(n + 4)} t_{(2)} \hit (g^0_{(2)} \cdots g^n_{(n + 2)}) 
t_{(3)}(\delta S^{-1}(f_{(1)})) \delta S^{-1}(f_{(2)}) \hitby \sigma \delta^{-1})(k) \otimes \\
&\qquad (f_{(3)}S^{-1}(g^0_{(1)} \cdots g^n_{(n + 1)}) \cdot  a^0) d(f_{(4)} S^{-1}(g^1_{(1)} \cdots g^n_{(n)})\cdot a^1) \cdots 
d(f_{(n + 3)} S^{-1}(g^n_{(1)}) \cdot a^n) \allowdisplaybreaks[2] \\ 
&= \hat{\F}(t_{(1)} \hit f_{(n + 3)} t_{(2)} \hit (g^0_{(2)} \cdots g^n_{(n + 2)}) 
t_{(3)} \hit (\delta S^{-1}(f_{(1)})) \hitby \sigma \delta^{-1})(k) \otimes \\
&\qquad  (f_{(2)}S^{-1}(g^0_{(1)} \cdots g^n_{(n + 1)}) \cdot  a^0) d(f_{(3)} S^{-1}(g^1_{(1)} \cdots g^n_{(n)})\cdot a^1) \cdots 
d(f_{(n + 2)} S^{-1}(g^n_{(1)}) \cdot a^n) \allowdisplaybreaks[2] \\ 
&= \hat{\psi}(\delta S^{-1}(\delta S^{-1}(f_{(1)}) \hitby \sigma) \hitby S(t_{(3)}) S^{-1}(g^0_{(2)} \cdots g^n_{(n + 2)}) \hitby S(t_{(2)}) 
S^{-1}(f_{(n + 3)}) \hitby S(t_{(1)}) \sigma \hit k \delta) \otimes \\
&\qquad (f_{(2)}S^{-1}(g^0_{(1)} \cdots g^n_{(n + 1)}) \cdot  a^0) 
d(f_{(3)} S^{-1}(g^1_{(1)} \cdots g^n_{(n)})\cdot a^1) \cdots d(f_{(n + 2)} S^{-1}(g^n_{(1)}) \cdot a^n) \allowdisplaybreaks[2] \\ 
&= \hat{\psi}((S^{-1}(\delta S^{-1}(f_{(1)}) \hitby \sigma)  S^{-1}(g^0_{(2)} \cdots g^n_{(n + 2)})   
S^{-1}(f_{(n + 3)}) (\sigma \hit k \delta^2) \hitby t_{(1)}) \hitby S(t_{(2)})) 
\otimes \\
&\qquad  (f_{(2)}S^{-1}(g^0_{(1)} \cdots g^n_{(n + 1)}) \cdot  a^0) d(f_{(3)} S^{-1}(g^1_{(1)} \cdots g^n_{(n)})\cdot a^1) \cdots 
d(f_{(n + 2)} S^{-1}(g^n_{(1)}) \cdot a^n),
\end{align*}
using $ \hat{\psi}(\delta h) = \hat{\psi}(h \delta) $ for all $ h \in \hat{H} $ in the last step. 
Since $ \hat{\psi} $ is right invariant this yields
\begin{align*}
&\quad \hat{\psi}(S^{-1}(\delta S^{-1}(f_{(1)}) \hitby \sigma)  S^{-1}(g^0_{(2)} \cdots g^n_{(n + 2)})   
S^{-1}(f_{(n + 3)}) (\sigma \hit k \delta^2) \hitby t_{(1)} \delta^{-2}(t_{(2)})) 
\otimes \\
&\qquad  (f_{(2)}S^{-1}(g^0_{(1)} \cdots g^n_{(n + 1)}) \cdot  a^0) d(f_{(3)} S^{-1}(g^1_{(1)} \cdots g^n_{(n)})\cdot a^1) \cdots 
d(f_{(n + 2)} S^{-1}(g^n_{(1)}) \cdot a^n) \allowdisplaybreaks[2] \\ 
&= \hat{\psi}(S^{-1}(\delta S^{-1}(f_{(1)}) \hitby \sigma)  S^{-1}(g^0_{(2)} \cdots g^n_{(n + 2)})   
S^{-1}(f_{(n + 3)}) \sigma \hit k_{(2)} \delta^2) t(k_{(1)}) 
\otimes \\
&\qquad  (f_{(2)}S^{-1}(g^0_{(1)} \cdots g^n_{(n + 1)}) \cdot  a^0) d(f_{(3)} S^{-1}(g^1_{(1)} \cdots g^n_{(n)})\cdot a^1) \cdots 
d(f_{(n + 2)} S^{-1}(g^n_{(1)}) \cdot a^n) \allowdisplaybreaks[2] \\ 
&= \hat{\F}(f_{(n + 3)} g^0_{(2)} \cdots g^n_{(n + 2)} \sigma \hit S(f_{(1)}))(k_{(2)}) t(k_{(1)}) 
\otimes \\
&\qquad  (f_{(2)}S^{-1}(g^0_{(1)} \cdots g^n_{(n + 1)}) \cdot  a^0) d(f_{(3)} S^{-1}(g^1_{(1)} \cdots g^n_{(n)})\cdot a^1) \cdots 
d(f_{(n + 2)} S^{-1}(g^n_{(1)}) \cdot a^n) \allowdisplaybreaks[2] \\ 
&= k_{(1)}(t) \tau_H(f \otimes (a^0 \rtimes g^0) 
d(a^1 \rtimes g^1) \cdots d(a^n \rtimes g^n))(k_{(2)}) \\
&= (t \bullet \tau_H(f \otimes (a^0 \rtimes g^0) d(a^1 \rtimes g^1) \cdots d(a^n \rtimes g^n)))(k),
\end{align*}
which proves that the map $ \tau_H $ is $ H $-linear. Moreover, for $ h \in \hat{H} $ we compute 
\begin{align*}
&\tau_H(h\bullet (f \otimes (a^0 \rtimes g^0) d(a^1 \rtimes g^1) \cdots d(a^n \rtimes g^n)))(k) \\
&= \tau_H(hf \otimes (a^0 \rtimes g^0) d(a^1 \rtimes g^1) \cdots d(a^n \rtimes g^n))(k) \\
&= \hat{\F}(h_{(n + 3)} f_{(n + 3)} g^0_{(2)} \cdots g^n_{(n + 2)} 
\sigma \hit S(f_{(1)}) \sigma \hit S(h_{(1)}))(k) \otimes \\
&\qquad \otimes (h_{(2)} f_{(2)} S^{-1}(g^0_{(1)} \cdots g^n_{(n + 1)}) \cdot  a^0) 
d(h_{(3)} f_{(3)} S^{-1}(g^1_{(1)} \cdots g^n_{(n)})\cdot a^1) \cdots \\
&\qquad \cdots d(h_{(n + 2)}f_{(n + 2)} S^{-1}(g^n_{(1)}) \cdot a^n) \allowdisplaybreaks[2] \\
&= \hat{\psi}(h_{(1)} \hitby \sigma^{-1} f_{(1)} \hitby \sigma^{-1} S^{-1}(g^0_{(2)} \cdots g^n_{(n + 2)})
S^{-1}(f_{(n + 3)}) S^{-1}(h_{(n + 3)}) \sigma \hit k \delta) \otimes \\
&\qquad \otimes (h_{(2)} f_{(2)} S^{-1}(g^0_{(1)} \cdots g^n_{(n + 1)}) \cdot  a^0) 
d(h_{(3)} f_{(3)} S^{-1}(g^1_{(1)} \cdots g^n_{(n)})\cdot a^1) \cdots \\
&\qquad \cdots d(h_{(n + 2)}f_{(n + 2)} S^{-1}(g^n_{(1)}) \cdot a^n), 
\end{align*}
and using the twisted trace property of $ \hat{\psi} $ we see that this is equal to
\begin{align*}
&\;\;\; \hat{\psi}(f_{(1)} \hitby \sigma^{-1} S^{-1}(g^0_{(2)} \cdots g^n_{(n + 2)})
S^{-1}(f_{(n + 3)}) S^{-1}(h_{(n + 3)}) \sigma \hit k \delta \sigma \hit (\delta^{-1} h_{(1)} \delta)) \otimes \\
&\qquad \otimes (h_{(2)} f_{(2)} S^{-1}(g^0_{(1)} \cdots g^n_{(n + 1)}) \cdot  a^0) 
d(h_{(3)} f_{(3)} S^{-1}(g^1_{(1)} \cdots g^n_{(n)})\cdot a^1) \cdots \\
&\qquad \cdots d(h_{(n + 2)}f_{(n + 2)} S^{-1}(g^n_{(1)}) \cdot a^n) \allowdisplaybreaks[2] \\
&= \hat{\psi}(f_{(1)} \hitby \sigma^{-1} S^{-1}(g^0_{(2)} \cdots g^n_{(n + 2)})
S^{-1}(f_{(n + 3)}) \sigma \hit S^{-1}(h_{(n + 4)}) \sigma \hit k \sigma \hit h_{(1)} \delta) \otimes \\
&\qquad \otimes (h_{(2)} f_{(2)} S^{-1}(g^0_{(1)} \cdots g^n_{(n + 1)}) \cdot  a^0) 
d(h_{(3)} f_{(3)} S^{-1}(g^1_{(1)} \cdots g^n_{(n)})\cdot a^1) \cdots \\
&\qquad \cdots d(h_{(n + 2)}f_{(n + 2)} S^{-1}(g^n_{(1)}) \cdot a^n) \sigma(h_{(n + 3)}) \allowdisplaybreaks[2] \\
&= \hat{\F}(f_{(n + 3)} g^0_{(2)} \cdots g^n_{(2)} \sigma \hit S(f_{(1)}))
(S^{-1}(h_{(4)}) k h_{(1)}) \otimes \\
&\qquad \otimes h_{(2)} \sigma(h_{(3)}) \cdot (f_{(2)} S^{-1}(g^0_{(1)} \cdots g^n_{(n + 1)}) \cdot  a^0) 
d(f_{(3)} S^{-1}(g^1_{(1)} \cdots g^n_{(n)})\cdot a^1) \cdots \\
&\qquad \cdots d(f_{(n + 2)} S^{-1}(g^n_{(1)}) \cdot a^n) \\
&= h \bullet \tau_H(f \otimes (a^0 \rtimes g^0) d(a^1 \rtimes g^1) \cdots d(a^n \rtimes g^n))(k).
\end{align*}
This shows that $ \tau_H $ is $ \hat{H} $-linear. \\
By construction, the map $ \tau_H $ is a chain map with respect to the operator $ d $. 
In order to show that $ \tau_H $ is a chain map with respect to the equivariant Hochschild boundary $ b_H $ we shall 
verify that $ \tau_H $ commutes with all simplicial boundary operators individually. 
Indeed, for $ 0 \leq j < n $ we have 
\begin{align*}
&\tau_H d_j^H(f \otimes (a^0 \rtimes g^0) d(a^1 \rtimes g^1) \cdots d(a^n \rtimes g^n)) \\
&= \tau_H(f \otimes (a^0 \rtimes g^0) d(a^1 \rtimes g^1) \cdots d(a^j g^j_{(1)} \cdot a^{j + 1} \rtimes g^j_{(2)} g^{j + 1}) \cdots d(a^n \rtimes g^n)) \\
&= \hat{\F}(f_{(n + 2)} g^0_{(2)} g^1_{(3)} \cdots g^j_{(j + 3)} g^{j + 1}_{(j + 2)} g^{j + 2}_{(j + 3)} \cdots 
g^n_{(n + 1)} \sigma \hit S(f_{(1)})) \otimes \\
&\qquad \otimes (f_{(2)}S^{-1}(g^0_{(1)} g^1_{(2)} \cdots g^j_{(j + 2)} g^{j + 1}_{(j + 1)} \cdots g^n_{(n)}) \cdot  a^0) 
d(f_{(3)}S^{-1}(g^1_{(1)} \cdots g^n_{(n - 1)}) \cdot a^1) \cdots \\
&\qquad \cdots d(f_{(j + 2)}S^{-1}(g^j_{(2)} g^{j + 1}_{(1)} \cdots g^n_{(n - j)}) \cdot  (a^j g^j_{(1)} \cdot a^{j + 1}))
\cdots d(f_{(n + 1)}S^{-1}(g^n_{(1)}) \cdot a^n) \allowdisplaybreaks[2]\\
&= \hat{\F}(f_{(n + 3)} g^0_{(2)} g^1_{(3)} \cdots g^j_{(j + 2)} g^{j + 1}_{(j + 3)} g^{j + 2}_{(j + 4)} \cdots g^n_{(n + 2)} 
\sigma \hit S(f_{(1)})) \otimes \\
&\qquad \otimes (f_{(2)}S^{-1}(g^0_{(1)} g^1_{(2)} \cdots g^n_{(n + 1)}) \cdot  a^0) 
d(f_{(3)}S^{-1}(g^1_{(1)} \cdots g^n_{(n)}) \cdot a^1) \cdots \\
&\qquad \cdots d(f_{(j + 2)}S^{-1}(g^j_{(1)} g^{j + 1}_{(2)} \cdots g^n_{(n - j + 1)}) \cdot  a^j 
f_{(j + 3)}S^{-1}(g^{j + 1}_{(1)}\cdots g^n_{(n - j)}) \cdot  a^{j + 1}) \cdots \\
&\qquad \cdots d(f_{(n + 2)}S^{-1}(g^n_{(1)}) \cdot  a^n) \\
&= d_j^{\hat{H}} \tau_H(f \otimes (a^0 \rtimes g^0) d(a^1 \rtimes g^1) \cdots d(a^n \rtimes g^n)).
\end{align*}
In the most interesting case $ j = n $ we calculate 
\begin{align*}
&\tau_H d^H_n (f \otimes (a^0 \rtimes g^0) d(a^1 \rtimes g^1) \cdots d(a^n \rtimes g^n))(k) \\
&= \tau_H(fS^{-1}(g^n_{(2)}) \hitby \sigma \otimes (a^n \rtimes g^n_{(1)})(a^0 \rtimes g^0) d(a^1 \rtimes g^1) \cdots d(a^{n - 1} \rtimes g^{n - 1}))(k) \\
&= \tau_H(fS^{-1}(g^n_{(3)}) \hitby \sigma \otimes 
(a^n g^n_{(1)} \cdot a^0 \rtimes g^n_{(2)} g^0) d(a^1 \rtimes g^1) \cdots d(a^{n - 1} \rtimes g^{n - 1}))(k) \\
&= \hat{\F}(f_{(n + 2)} S^{-1}(g^n_{(4)}) g^n_{(3)} g^0_{(2)} g^1_{(3)} \cdots g^{n - 1}_{(n + 1)} 
\sigma \hit S(f_{(1)} S^{-1}(g^n_{(n + 5)}) \hitby \sigma))(k) \otimes \allowdisplaybreaks[2] \\
&\qquad f_{(2)} S^{-1}(g^n_{(n + 4)}) S^{-1}(g^n_{(2)} g^0_{(1)} \cdots g^{n - 1}_{(n)}) \cdot (a^n g^n_{(1)} \cdot a^0) \\
&\qquad d(f_{(3)} S^{-1}(g^n_{(n + 3)}) S^{-1}(g^1_{(1)} \cdots g^{n - 1}_{(n - 1)}) \cdot a^1) \cdots 
d(f_{(n + 1)} S^{-1}(g^n_{(5)}) S^{-1}(g^{n - 1}_{(1)}) \cdot a^{n - 1}) \allowdisplaybreaks[2] \\
&= \hat{\F}(f_{(n + 2)} g^0_{(2)} g^1_{(3)} \cdots g^{n - 1}_{(n + 1)} g^n_{(n + 3)} \sigma \hit S(f_{(1)}))(k) \otimes \\
&\qquad f_{(2)} S^{-1}(g^n_{(2)} g^0_{(1)} \cdots g^{n - 1}_{(n)} g^n_{(n + 2)}) \cdot (a^n g^n_{(1)} \cdot a^0) \\
&\qquad d(f_{(3)} S^{-1}(g^1_{(1)} \cdots g^{n - 1}_{(n - 1)} g^n_{(n + 1)}) \cdot a^1) \cdots 
d(f_{(n + 1)} S^{-1}(g^{n - 1}_{(1)} g^n_{(3)}) \cdot a^{n - 1}) \allowdisplaybreaks[2] \\
&= \hat{\F}(f_{(n + 2)} g^0_{(2)} g^1_{(3)} \cdots g^{n - 1}_{(n + 1)} g^n_{(n + 2)} \sigma \hit S(f_{(1)}))(k) \otimes \\
&\qquad f_{(2)} S^{-1}(g^0_{(1)} \cdots g^{n - 1}_{(n)} g^n_{(n + 1)}) \cdot ((S^{-1}(g^n_{(1)}) \cdot a^n) a^0) \\
&\qquad d(f_{(3)} S^{-1}(g^1_{(1)} \cdots g^{n - 1}_{(n - 1)} g^n_{(n)}) \cdot a^1) \cdots 
d(f_{(n + 1)} S^{-1}(g^{n - 1}_{(1)} g^n_{(2)}) \cdot a^{n - 1}) \\
&= \hat{\F}(f_{(n + 3)} g^0_{(3)} g^1_{(4)} \cdots g^{n - 1}_{(n + 2)} g^n_{(n + 3)} \sigma \hit S(f_{(1)}))(k) \otimes \\
&\qquad (f_{(2)} S^{-1}(g^0_{(2)} \cdots g^{n - 1}_{(n + 1)} g^n_{(n + 2)}) \cdot (S^{-1}(g^n_{(1)}) \cdot a^n) 
(f_{(3)} S^{-1}(g^0_{(1)} \cdots g^{n - 1}_{(n)} g^n_{(n + 1)}) \cdot a^0) \\
&\qquad d(f_{(4)} S^{-1}(g^1_{(1)} \cdots g^{n - 1}_{(n - 1)} g^n_{(n)}) \cdot a^1) \cdots 
d(f_{(n + 2)} S^{-1}(g^{n - 1}_{(1)} g^n_{(2)}) \cdot a^{n - 1}) \allowdisplaybreaks[2] \\
&= \hat{\psi}(S^{-1}(\delta^{-1} S(k) \hitby \sigma^{-1} f_{(n + 5)} g^0_{(3)} \cdots g^n_{(n + 3)} \sigma \hit S(f_{(1)}))) \otimes \\
&\qquad (S^{-1}(f_{(n + 4)} g^0_{(2)} \cdots g^n_{(n + 2)} S(f_{(2)})) f_{(n + 3)} S^{-1}(g^n_{(1)}) \cdot a^n) 
(f_{(3)} S^{-1}(g^0_{(1)} \cdots g^{n - 1}_{(n)} g^n_{(n + 1)}) \cdot a^0)  \\
&\qquad d(f_{(4)} S^{-1}(g^1_{(1)} \cdots g^{n - 1}_{(n - 1)} g^n_{(n)}) \cdot a^1) \cdots 
d(f_{(n + 2)} S^{-1}(g^{n - 1}_{(1)} g^n_{(2)}) \cdot a^{n - 1}) \allowdisplaybreaks[2] \\
&= \hat{\psi}(S^{-1}(\delta^{-1} S(k_{(1)}) f_{(n + 5)} g^0_{(3)} \cdots g^n_{(n + 3)} \sigma \hit S(f_{(1)}))) \otimes \\
&\qquad (S^{-1}(\sigma \hit S^2(k_{(3)}) \delta \delta^{-1} S(k_{(2)}) f_{(n + 4)} g^0_{(2)} \cdots 
g^n_{(n + 2)} S(f_{(2)})) f_{(n + 3)} S^{-1}(g^n_{(1)}) \cdot a^n) \\
&\qquad (f_{(3)} S^{-1}(g^0_{(1)} \cdots g^{n - 1}_{(n)} g^n_{(n + 1)}) \cdot a^0)
d(f_{(4)} S^{-1}(g^1_{(1)} \cdots g^{n - 1}_{(n - 1)} g^n_{(n)}) \cdot a^1) \cdots \\
&\qquad \cdots d(f_{(n + 2)} S^{-1}(g^{n - 1}_{(1)} g^n_{(2)}) \cdot a^{n - 1}). 
\end{align*}
Using that $ \hat{\psi} $ is right invariant this becomes
\begin{align*}
&\qquad \hat{\psi}(S^{-1}(\delta^{-1} S(k_{(1)}) f_{(n + 3)} g^0_{(2)} \cdots g^n_{(n + 2)} \sigma \hit S(f_{(1)}))) \otimes \\
&\qquad (S^{-1}(\sigma \hit S^2(k_{(2)}) \delta) 
f_{(n + 2)} S^{-1}(g^n_{(1)}) \cdot a^n) (f_{(2)} S^{-1}(g^0_{(1)} \cdots g^{n}_{(n + 1)}) \cdot a^0) \\
&\qquad d(f_{(3)} S^{-1}(g^1_{(1)} \cdots g^{n - 1}_{(n - 1)} g^n_{(n)}) \cdot a^1) \cdots 
d(f_{(n + 1)} S^{-1}(g^{n - 1}_{(1)} g^n_{(2)}) \cdot a^{n - 1}) \allowdisplaybreaks[2]\\
&= \hat{\psi}(S^{-1}(\delta^{-1} S(k_{(1)}) f_{(n + 3)} g^0_{(2)} \cdots g^n_{(n + 2)} \sigma \hit S(f_{(1)}))) \otimes \\
&\qquad (S^{-1}(\delta k_{(2)} \hitby \sigma) 
f_{(n + 2)} S^{-1}(g^n_{(1)}) \cdot a^n) (f_{(2)} S^{-1}(g^0_{(1)} \cdots g^{n}_{(n + 1)}) \cdot a^0) \\
&\qquad d(f_{(3)} S^{-1}(g^1_{(1)} \cdots g^{n - 1}_{(n - 1)} g^n_{(n)}) \cdot a^1) \cdots 
d(f_{(n + 1)} S^{-1}(g^{n - 1}_{(1)} g^n_{(2)}) \cdot a^{n - 1}) \allowdisplaybreaks[2]\\
&= \hat{\psi}(S^{-1}(\delta^{-1} S(k_{(1)}) \hitby \sigma^{-1} f_{(n + 3)} g^0_{(2)} \cdots g^n_{(n + 2)} \sigma \hit S(f_{(1)}))) \otimes \\
&\qquad (S^{-1}(\delta k_{(2)}) 
f_{(n + 2)} S^{-1}(g^n_{(1)}) \cdot a^n) (f_{(2)} S^{-1}(g^0_{(1)} \cdots g^{n}_{(n + 1)}) \cdot a^0) \\
&\qquad d(f_{(3)} S^{-1}(g^1_{(1)} \cdots g^{n - 1}_{(n - 1)} g^n_{(n)}) \cdot a^1) \cdots 
d(f_{(n + 1)} S^{-1}(g^{n - 1}_{(1)} g^n_{(2)}) \cdot a^{n - 1}) \\
&= d_n(\hat{\F}(f_{(n + 3)} g^0_{(2)} g^1_{(3)} \cdots g^n_{(n + 2)} \sigma \hit S(f_{(1)})) \otimes 
f_{(2)} S^{-1}(g^0_{(1)} g^1_{(2)} \cdots g^n_{(n + 1)}) \cdot a^0 \\
&\qquad d(f_{(3)} S^{-1}(g^1_{(1)} \cdots g^n_{(n)}) \cdot a^1) \cdots 
d(f_{(n + 2)} S^{-1}(g^n_{(1)}) \cdot a^n))(k) \\
&= d_n^{\hat{H}} \tau_H(f \otimes (a^0 \rtimes g^0) d(a^1 \rtimes g^1) \cdots d(a^n \rtimes g^n))(k).
\end{align*}
Since the equivariant Hochschild and Connes boundary operators on both sides are constructed out of $ d $ and the simplicial boundaries this finishes 
the proof. \qed \\
We need another lemma.
\begin{lemma} \label{lambdaflemma} 
Let $ \F: H \rightarrow \hat{H} $ denote the twisted Fourier transform for the dual. 
Then for $ t \in H $ and $ f \in \hat{H} $ the following relations hold. 
\begin{bnum}
\item[a)] $ \hat{\G}_l(\sigma \hit \F(t) \delta) = t $.  
\item[b)] $ f \bullet \F(t) = f(\hat{\psi}_{(2)} \sigma) \F(t_{(1)}) \psi(\hat{\psi}_{(1)} S^{-1}(t_{(2)})) $.
\end{bnum}
\end{lemma}
\proof a) For $ h \in \hat{H} $ we calculate
\begin{align*}
(h \sigma \hit \F(t) \delta)(r) &= h(r_{(1)}) \F(t)(r_{(2)} \sigma) \delta(r_{(3)}) \\
&= h(r_{(1)}) \psi(S^{-1}(t) \delta \hit r_{(2)} \sigma^2) \delta(r_{(3)}) \\
&= h(t_{(3)} S^{-1}(t_{(2)}) r_{(1)} \sigma^2 \sigma^{-2}) \psi(S^{-1}(t_{(1)}) r_{(2)} \sigma^2) \delta^2(r_{(3)}) \\
&= h(t_{(2)} \sigma^2 \sigma^{-2}) \psi(S^{-1}(t_{(1)}) r_{(1)} \sigma^2) \delta^2(r_{(2)}) \\
&= h(t_{(4)}) \psi(S^{-1}(t_{(3)}) r_{(1)} \sigma^2) \delta^2(S^{-2}(t_{(1)}) S^{-1}(t_{(2)}) r_{(2)} \sigma^2 \sigma^{-2}) \\
&= h(t_{(3)}) \psi(S^{-1}(t_{(2)}) r \sigma^2) \delta^2(S^{-2}(t_{(1)}) \sigma^{-2}) \\
&= h(t_{(3)}) \psi(r \sigma \delta \hit S^{-1}(t_{(2)}) \hitby \delta \sigma) \delta^2(t_{(1)}), 
\end{align*}
and hence we get
\begin{align*}
\hat{\G}_l(\sigma \hit \F(t) \delta)(h) &= \hat{\psi}(h \sigma \hit \F(t) \delta) \\
&= h(t_{(3)}) \epsilon(\sigma \delta \hit S^{-1}(t_{(2)}) \hitby \delta \sigma) \delta^2(t_{(1)}) \\
&= h(t).
\end{align*}
Here we have used that $ \psi $ is a left invariant integral for $ H^{\cop} $, the bornological quantum group $ H $ equipped with the 
opposite comultiplication, and we note that $ \hat{\psi} $ is a right invariant integral for the dual of $ H^{\cop} $. \\
b) We compute 
\begin{align*}
(f \bullet \F(t))(r) & = f(r_{(1)}) \psi(S^{-1}(t) \delta \hit r_{(2)} \sigma) \\
&= \delta(r_{(3)}) f(t_{(3)} S^{-1}(t_{(2)}) r_{(1)} \sigma \sigma^{-1}) \psi(S^{-1}(t_{(1)}) r_{(2)} \sigma) \\
&= \delta(r_{(2)}) f(t_{(2)} \sigma^2 \sigma^{-1}) \psi(S^{-1}(t_{(1)}) r_{(1)} \sigma) \\
&= f(t_{(2)} \sigma) \F(t_{(1)})(r)  
\end{align*}
and 
\begin{align*}
f(t_{(2)} \sigma) \F(t_{(1)}) 
&= f(\hat{\psi}_{(2)} t_{(2)} \sigma) \F(t_{(1)}) \psi(\hat{\psi}_{(1)}) \\
&= f(\hat{\psi}_{(2)} S^{-1}(t_{(3)}) t_{(2)} \sigma) \F(t_{(1)}) \psi(\hat{\psi}_{(1)} S^{-1}(t_{(4)})) \\
&= f(\hat{\psi}_{(2)} \sigma) \F(t_{(1)}) \psi(\hat{\psi}_{(1)} S^{-1}(t_{(2)})). 
\end{align*}
This proves the claim. \qed \\
The following main technical result shows that the duality map introduced above provides a 
natural factorisation of the trace map $ \tr_A: \Omega_H(A \cotimes \mathcal{K}_H) \rightarrow \Omega_H(A) $, defined by 
\begin{align*}
\tr_A(t \otimes &(a^0 \otimes r^0 \otimes g^0) d(a^1 \otimes r^1 \otimes g^1) \cdots d(a^n \otimes r^n \otimes g^n)) \\
&= t_{(2)} \otimes g^n(t_{(1)} \cdot r^0) g^0(r^1) \cdots g^{n - 1}(r^n) a^0 da^1 \cdots da^n. 
\end{align*}
We remark that this map plays an important role in the proof of stability \cite{Voigtechq}. 
\begin{theorem} \label{technical}
Let $ A $ be any $ H $-algebra. For $ n = 0 $ and $ n = 1 $ the map 
$$ 
\tau_\ihd \tau_H: \Omega_H^n(A \rtimes H \rtimes \hat{H}) \rightarrow \Omega^n_H(A) 
$$ 
is equal to the composition of the natural automorphism $ T $ with the Takesaki-Takai isomorphism $ \Omega_H^n(A \rtimes H \rtimes \hat{H}) \cong 
\Omega_H^n(A \cotimes \mathcal{K}_H) $ and the trace map $ \tr_A: \Omega_H^n(A \cotimes \mathcal{K}_H) \rightarrow \Omega^n_H(A) $. 
\end{theorem}
\proof Using proposition \ref{duality1} we compute 
\begin{align*}
&\tau_\ihd \tau_H(f \otimes a^0 \rtimes r^0 \rtimes g^0) = \tau_\ihd(\hat{\F}(f_{(3)} g^0_{(2)} \sigma \hit S(f_{(1)})) \otimes 
(f_{(2)}S^{-1}(g^0_{(1)}))(r^0_{(2)}) a^0 \rtimes r^0_{(1)}) \\
&= \tau_\ihd(f \bullet (\hat{\F}(g^0_{(2)}) \otimes S^{-1}(g^0_{(1)})(r^0_{(2)}) a^0 \rtimes r^0_{(1)})) \\
&= f \bullet \tau_\ihd(\hat{\psi}_{(1)}(\delta S^{-1}(g^0_{(2)})) \hat{\psi}_{(2)} \sigma \otimes S^{-1}(g^0_{(1)})(r^0_{(2)}) a^0 \rtimes r^0_{(1)})\\
&= f \bullet \F(\hat{\psi}_{(4)} \sigma r^0_{(2)} \sigma^{-1} \delta \hit S(\hat{\psi}_{(2)})) 
\otimes \hat{\psi}_{(1)}(\delta S^{-1}(g^0_{(2)})) S^{-1}(g^0_{(1)})(r^0_{(3)}) (\hat{\psi}_{(3)} \sigma S^{-1}(r^0_{(1)})) \cdot a^0 \\
&= S^{-1}(g^0)(\hat{\psi}_{(2)} r^0_{(3)}) \delta(\hat{\psi}_{(1)}) f \bullet \F(\hat{\psi}_{(5)} \sigma r^0_{(2)} \sigma^{-1} \delta \hit S(\hat{\psi}_{(3)})) 
\otimes (\hat{\psi}_{(4)} \sigma S^{-1}(r^0_{(1)})) \cdot a^0. 
\end{align*}
Applying the operator $ T^{-1} $ and using lemma \ref{lambdaflemma} a) we obtain
\begin{align*}
T^{-1} &\tau_\ihd \tau_H (f \otimes a^0 \rtimes r^0 \rtimes g^0) \\
&= S^{-1}(g^0)(\hat{\psi}_{(2)} r^0_{(4)}) \delta(\hat{\psi}_{(1)}) f \bullet \F(\hat{\psi}_{(7)} \sigma r^0_{(3)} \sigma^{-1} \delta \hit S(\hat{\psi}_{(3)})) \\
&\qquad \otimes (\hat{\psi}_{(6)} \sigma r^0_{(2)} \sigma^{-1} S(\hat{\psi}_{(4)}) \hat{\psi}_{(5)} \sigma S^{-1}(r^0_{(1)})) \cdot a^0 \\
&= S^{-1}(g^0)(\hat{\psi}_{(2)} r^0_{(2)}) \delta(\hat{\psi}_{(1)}) 
f \bullet \F(\hat{\psi}_{(5)} \sigma r^0_{(1)} \sigma^{-1} \delta \hit S(\hat{\psi}_{(3)})) (\hat{\psi}_{(4)} \sigma) \cdot a^0 \\
&= S^{-1}(g^0)(\hat{\psi}_{(1)} r^0_{(2)}) \delta(\hat{\psi}_{(5)}) 
f \bullet \F(\hat{\psi}_{(4)} \sigma r^0_{(1)} \sigma^{-1} \delta \hit S(\hat{\psi}_{(2)})) \otimes (\hat{\psi}_{(3)} \sigma) \cdot a^0 \\
&= S^{-1}(g^0)(\hat{\psi}_{(1)} \epsilon(r^0_{(6)}) r^0_{(7)}) \delta(\hat{\psi}_{(5)} \epsilon(r^0_{(2)})) \\
&\qquad f \bullet \F(\hat{\psi}_{(4)} \epsilon(r^0_{(3)}) \sigma r^0_{(1)} \sigma^{-1} \delta \hit S(\hat{\psi}_{(2)} \epsilon(r^0_{(5)}))) 
\otimes (\hat{\psi}_{(3)} \epsilon(r^0_{(4)}) \sigma) \cdot a^0, 
\end{align*}
where we use $ \hat{\psi}(\delta h) = \hat{\psi}(h \delta) $ for all $ h \in \hat{H} $ in the last equality. 
Using that $ \hat{\psi} $ is a right invariant integral we obtain 
\begin{align*}
&T^{-1} \tau_\ihd \tau_H(f \otimes a^0 \rtimes r^0 \rtimes g^0) = 
S^{-1}(g^0)(\hat{\psi}_{(1)} S(r^0_{(5)}) r^0_{(6)}) \delta(\hat{\psi}_{(5)}) \\
&\qquad f \bullet \F(\hat{\psi}_{(4)} \delta \hit S(r^0_{(2)}) \sigma r^0_{(1)} \sigma^{-1} \delta \hit S(\hat{\psi}_{(2)} S(r^0_{(4)}))) 
\otimes (\hat{\psi}_{(3)} S(r^0_{(3)})\sigma) \cdot a^0 \\
&= S^{-1}(g^0)(\hat{\psi}_{(1)}) \delta(\hat{\psi}_{(5)}) 
f \bullet \F(\hat{\psi}_{(4)} \sigma S^{-1}(r^0_{(2)}) \hitby \delta  r^0_{(1)} \sigma^{-1} \delta \hit S(\hat{\psi}_{(2)} S(r^0_{(4)}))) \\
&\qquad \otimes (\hat{\psi}_{(3)} S(r^0_{(3)}) \sigma) \cdot a^0 \\
&= S^{-1}(g^0)(\hat{\psi}_{(1)}) \delta(\hat{\psi}_{(5)}) 
f \bullet \F(\hat{\psi}_{(4)} \delta^{-1}(r^0_{(1)}) \delta \hit S(\hat{\psi}_{(2)} S(r^0_{(3)}))) \\
&\qquad \otimes (\hat{\psi}_{(3)} S(r^0_{(2)}) \sigma) \cdot a^0.
\end{align*}
In the following computation we have to diagonalise two right invariant integrals $ \hat{\psi} $. In order to distinguish 
them we will use square brackets to denote Sweedler indices for the first one. 
Using lemma \ref{lambdaflemma} a) and the Takesaki-Takai isomorphism from proposition \ref{TTduality} we calculate 
\begin{align*}
&\tr_A(f \otimes a^0 \rtimes r^0 \rtimes g^0) = \tr_A(f(\hat{\psi}_{[2]} \sigma) \hat{\psi}_{[3]}(\delta) \F(\hat{\psi}_{[1]}) \otimes 
\gamma_A(a^0 \rtimes r^0 \rtimes g^0)) \\
&= S^{-1}(g^0)(\hat{\psi}_{(1)}) f(\hat{\psi}_{[3]} \sigma) \delta(\hat{\psi}_{[4]}) \F(\hat{\psi}_{[2]}) \\
&\qquad \psi(\hat{\psi}_{[1]} \hat{\psi}_{(2)} S(r^0_{(2)}) S^{-1}(\hat{\psi}_{(4)})) 
\otimes (\hat{\psi}_{(3)} \delta \hit S(r^0_{(1)}) \sigma) \cdot a^0 \allowdisplaybreaks[2] \\
&= S^{-1}(g^0)(\hat{\psi}_{(1)}) f(\hat{\psi}_{[3])} \sigma) \delta(\hat{\psi}_{[4]}) \F(\hat{\psi}_{[2]}) \\
&\qquad \psi(\hat{\psi}_{[1]} \hat{\psi}_{(2)} \delta \hit S(r^0_{(2)}) S^{-1}(\hat{\psi}_{(4)})) 
\otimes (\hat{\psi}_{(3)} \sigma S^{-1}(r^0_{(1)})) \cdot a^0. 
\end{align*}
Using $ \hat{\psi}(\delta h) = \hat{\psi}(h \delta) $ for all $ h \in \hat{H} $ this gives
\begin{align*}
&\quad S^{-1}(g^0)(\hat{\psi}_{(1)}) f(\hat{\psi}_{[3]} \sigma)
\F(\hat{\psi}_{[2]}) \\
&\qquad \psi(\hat{\psi}_{[1]} \hitby \delta \hat{\psi}_{(2)} \delta \hit S(r^0_{(2)}) S^{-1}(\hat{\psi}_{(4)})) 
\otimes (\hat{\psi}_{(3)} \sigma S^{-1}(r^0_{(1)})) \cdot a^0 \allowdisplaybreaks[2] \\
&= S^{-1}(g^0)(\hat{\psi}_{(1)}) f(\hat{\psi}_{[3]} \sigma)
\F(\hat{\psi}_{[2]}) \\
&\qquad \psi(\hat{\psi}_{[1]} \hat{\psi}_{(2)} \hitby \delta^{-1} \sigma S^{-1}(r^0_{(2)}) \sigma^{-1} S^{-1}(\hat{\psi}_{(4)})\hitby \delta^{-1}) 
\otimes (\hat{\psi}_{(3)} \sigma S^{-1}(r^0_{(1)})) \cdot a^0 
\end{align*}
by invariance of $ \psi $ and $ \delta(\sigma) = 1 $. This equals
\begin{align*}
&\quad S^{-1}(g^0)(\hat{\psi}_{(1)}) f(\hat{\psi}_{[3]} \sigma)
\F(\hat{\psi}_{[2]} \hat{\psi}_{(3)} \sigma S^{-1}(r^0_{(3)}) \sigma^{-1} S^{-1}(\hat{\psi}_{(7)}) \hat{\psi}_{(6)} \sigma r^0_{(2)} \sigma^{-1} S(\hat{\psi}_{(4)})) \\
&\qquad \psi(\hat{\psi}_{[1]} \hat{\psi}_{(2)} \hitby \delta^{-1} \sigma S^{-1}(r^0_{(4)}) \sigma^{-1} S^{-1}(\hat{\psi}_{(8)}) \hitby \delta^{-1}) 
\otimes (\hat{\psi}_{(5)} \sigma S^{-1}(r^0_{(1)})) \cdot a^0 \allowdisplaybreaks[2] \\
&= S^{-1}(g^0)(\hat{\psi}_{(1)}) f(\hat{\psi}_{[2]} \sigma) \\
&\qquad \F(\hat{\psi}_{(5)} \sigma r^0_{(2)} \sigma^{-1} S(\hat{\psi}_{(3)})) 
\psi(\hat{\psi}_{[1]} \hat{\psi}_{(2)} \hitby \delta^{-1} \sigma S^{-1}(r^0_{(3)}) \sigma^{-1} S^{-1}(\hat{\psi}_{(6)})\hitby \delta^{-1}) \\
&\qquad \otimes (\hat{\psi}_{(4)} \sigma S^{-1}(r^0_{(1)})) \cdot a^0
\end{align*}
by right invariance of $ \psi $, which simplifies to
\begin{align*}
&S^{-1}(g^0)(\hat{\psi}_{(1)}) f(\hat{\psi}_{[2]} \sigma)\\
&\quad \F(\hat{\psi}_{(5)} \sigma r^0_{(2)} \sigma^{-1} S(\hat{\psi}_{(3)})) 
\psi(\hat{\psi}_{[1]} S^{-1}(\delta \hit \hat{\psi}_{(6)} \sigma r^0_{(3)} \sigma^{-1} \delta \hit S(\hat{\psi}_{(2)}))) \\
&\qquad \otimes (\hat{\psi}_{(4)} \sigma S^{-1}(r^0_{(1)})) \cdot a^0 \allowdisplaybreaks[2] \\
&= S^{-1}(g^0)(\hat{\psi}_{(1)}) f \bullet \F(\delta \hit \hat{\psi}_{(4)} \sigma r^0_{(2)} \sigma^{-1} \delta \hit S (\hat{\psi}_{(2)})) 
\otimes (\hat{\psi}_{(3)} \sigma S^{-1}(r^0_{(1)})) \cdot a^0 
\end{align*}
by lemma \ref{lambdaflemma} b). Hence we arrive at
\begin{align*} 
&\quad S^{-1}(g^0)(\hat{\psi}_{(1)}) f \bullet \F(\delta \hit \hat{\psi}_{(4)} \delta \hit S^2(r^0_{(2)}) \hitby \delta^{-1} \delta \hit 
S(\hat{\psi}_{(2)})) \\
&\qquad \otimes (\hat{\psi}_{(3)} \delta \hit S(r^0_{(1)}) \hitby \delta^{-1} \sigma) \cdot a^0 \\
&= S^{-1}(g^0)(\hat{\psi}_{(1)}) f \bullet \F(\delta \hit \hat{\psi}_{(4)} \delta \hit S^2(r^0_{(2)}) \delta \hit S(\hat{\psi}_{(2)})) \\
&\qquad \otimes (\hat{\psi}_{(3)} \delta \hit S(r^0_{(1)}) \sigma) \cdot a^0 \\
&= S^{-1}(g^0)(\hat{\psi}_{(1)}) \delta(\hat{\psi}_{(5)}) 
f \bullet \F(\hat{\psi}_{(4)} \delta^{-1}(r^0_{(1)}) \delta \hit S(\hat{\psi}_{(2)} S(r^0_{(3)}))) \\
&\qquad \otimes (\hat{\psi}_{(3)} S(r^0_{(2)}) \sigma) \cdot a^0
\end{align*}
which shows that $ T^{-1} \tau_\ihd \tau_H $ and $ \tr_A $ agree on equivariant differential forms of degree zero. \\
For equivariant differential forms of degree one we compute 
\begin{align*}
&\tau_\ihd \tau_H (f \otimes (a^0 \rtimes r^0 \rtimes g^0)d(a^1 \rtimes r^1 \rtimes g^1)) = 
\tau_\ihd(\hat{\F}(f_{(4)} g^0_{(2)} g^1_{(3)} \sigma \hit S(f_{(1)})) \otimes \\
&\qquad (f_{(2)}S^{-1}(g^0_{(1)} g^1_{(2)}))(r^0_{(2)}) (a^0 \rtimes r^0_{(1)}) d(f_{(3)}S^{-1}(g^1_{(1)}))(r^1_{(2)}) (a^1 \rtimes r^1_{(1)})) \\
&= \tau_\ihd(f \bullet (\hat{\F}(g^0_{(2)}g^1_{(3)}) \otimes 
S^{-1}(g^0_{(1)} g^1_{(2)})(r^0_{(2)}) (a^0 \rtimes r^0_{(1)}) dS^{-1}(g^1_{(1)})(r^1_{(2)}) (a^1 \rtimes r^1_{(1)}))) \\
&= f \bullet \tau_\ihd(\hat{\F}(g^0_{(2)}g^1_{(3)}) \otimes 
S^{-1}(g^0_{(1)} g^1_{(2)})(r^0_{(2)}) (a^0 \rtimes r^0_{(1)}) dS^{-1}(g^1_{(1)})(r^1_{(2)}) (a^1 \rtimes r^1_{(1)})) \\
&= f \bullet \tau_\ihd(\hat{\psi}_{(1)}(\delta S^{-1}(g^0_{(2)} g^1_{(3)})) \hat{\psi}_{(2)} \sigma \otimes 
S^{-1}(g^0_{(1)} g^1_{(2)})(r^0_{(2)}) (a^0 \rtimes r^0_{(1)}) dS^{-1}(g^1_{(1)})(r^1_{(2)}) (a^1 \rtimes r^1_{(1)})) \\
&= f \bullet \F(\hat{\psi}_{(5)} \sigma r^0_{(2)} r^1_{(3)} \sigma^{-1} \delta \hit S(\hat{\psi}_{(2)})) 
\otimes \hat{\psi}_{(1)}(\delta S^{-1}(g^0_{(2)} g^1_{(3)})) \\
&\qquad S^{-1}(g^0_{(1)}g^1_{(2)})(r^0_{(3)}) S^{-1}(g^1_{(1)})(r^1_{(4)}) 
(\hat{\psi}_{(3)} \sigma S^{-1}(r^0_{(1)} r^1_{(2)}) \cdot a^0) 
d(\hat{\psi}_{(4)} \sigma S^{-1}(r^1_{(1)}) \cdot a^1) \\
&= S^{-1}(g^0)(\hat{\psi}_{(3)} r^0_{(4)}) S^{-1}(g^1)(\hat{\psi}_{(2)} r^0_{(3)} r^1_{(4)}) \delta(\hat{\psi}_{(1)})
f \bullet \F(\hat{\psi}_{(7)} \sigma r^0_{(2)} r^1_{(3)} \sigma^{-1} \delta \hit S(\hat{\psi}_{(4)})) \\
&\qquad \otimes (\hat{\psi}_{(5)} \sigma S^{-1}(r^0_{(1)} r^1_{(2)}) \cdot a^0) d(\hat{\psi}_{(6)} \sigma S^{-1}(r^1_{(1)}) \cdot a^1).
\end{align*}
Applying the operator $ T^{-1} $ we obtain using lemma \ref{lambdaflemma} a)
\begin{align*}
&T^{-1} \tau_\ihd \tau_H(f \otimes (a^0 \rtimes r^0 \rtimes g^0) d(a^1 \rtimes r^1 \rtimes g^1)) \\
&= S^{-1}(g^0)(\hat{\psi}_{(3)} r^0_{(6)}) S^{-1}(g^1)(\hat{\psi}_{(2)} r^0_{(5)} r^1_{(6)}) \delta(\hat{\psi}_{(1)}) \\
&\qquad f \bullet  \F(\hat{\psi}_{(11)} \sigma r^0_{(4)} r^1_{(5)} \sigma^{-1} \delta \hit S(\hat{\psi}_{(4)})) 
\otimes (\hat{\psi}_{(9)} \sigma r^0_{(2)} r^1_{(3)} \sigma^{-1} S(\hat{\psi}_{(6)}) \hat{\psi}_{(7)} \sigma S^{-1}(r^0_{(1)} r^1_{(2)}) \cdot a^0) \\
&\qquad d(\hat{\psi}_{(10)} \sigma r^0_{(3)} r^1_{(4)} \sigma^{-1} S(\hat{\psi}_{(5)}) \hat{\psi}_{(8)} \sigma S^{-1}(r^1_{(1)}) \cdot a^1) \\
&= S^{-1}(g^0)(\hat{\psi}_{(3)} r^0_{(4)}) S^{-1}(g^1)(\hat{\psi}_{(2)} r^0_{(3)} r^1_{(2)}) \delta(\hat{\psi}_{(1)}) \\
&\qquad f \bullet \F(\hat{\psi}_{(7)} \sigma r^0_{(2)} r^1_{(1)} \sigma^{-1} \delta \hit S(\hat{\psi}_{(4)})) 
\otimes (\hat{\psi}_{(5)} \sigma \cdot a^0) d(\hat{\psi}_{(6)} \sigma r^0_{(1)} \cdot a^1) \\
&= S^{-1}(g^0)(\hat{\psi}_{(2)} r^0_{(4)}) S^{-1}(g^1)(\hat{\psi}_{(1)} r^0_{(3)} r^1_{(2)}) \delta(\hat{\psi}_{(7)}) \\
&\qquad f \bullet \F(\hat{\psi}_{(6)} \sigma r^0_{(2)} r^1_{(1)} \sigma^{-1} \delta \hit S(\hat{\psi}_{(3)})) 
\otimes (\hat{\psi}_{(4)} \sigma \cdot a^0) d(\hat{\psi}_{(5)} \sigma r^0_{(1)} \cdot a^1),
\end{align*}
and using that $ \hat{\psi} $ is a right invariant integral yields 
\begin{align*}
T^{-1} &\tau_\ihd \tau_H(f \otimes (a^0 \rtimes r^0 \rtimes g^0) d(a^1 \rtimes r^1 \rtimes g^1)) = 
S^{-1}(g^0)(\hat{\psi}_{(2)} S(r^0_{(7)} r^1_{(6)}) r^0_{(10)}) \\
&\qquad S^{-1}(g^1)(\hat{\psi}_{(1)} S(r^0_{(8)} r^1_{(7)}) r^0_{(9)} r^1_{(8)}) \delta(\hat{\psi}_{(7)}) \\
&\qquad f \bullet \F(\hat{\psi}_{(6)} \delta \hit S(r^0_{(3)} r^1_{(2)}) \sigma r^0_{(2)} r^1_{(1)} \sigma^{-1} 
\delta \hit S(\hat{\psi}_{(3)} S(r^0_{(6)} r^1_{(5)}))) \\
&\qquad \otimes (\hat{\psi}_{(4)} S(r^0_{(5)} r^1_{(4)}) \sigma \cdot a^0) d(\hat{\psi}_{(5)} S(r^0_{(4)} r^1_{(3)}) \sigma r^0_{(1)} \cdot a^1) \\
&= S^{-1}(g^0)(\hat{\psi}_{(2)} S(r^1_{(5)})) S^{-1}(g^1)(\hat{\psi}_{(1)}) \delta(\hat{\psi}_{(7)}) \\
&\qquad f \bullet \F(\hat{\psi}_{(6)} \delta^{-1}(r^0_{(2)} r^1_{(1)}) \delta \hit S(\hat{\psi}_{(3)} S(r^0_{(5)} r^1_{(4)}))) \\
&\qquad \otimes (\hat{\psi}_{(4)} S(r^0_{(4)} r^1_{(3)}) \sigma \cdot a^0) d(\hat{\psi}_{(5)} S(r^0_{(3)} r^1_{(2)}) \sigma r^0_{(1)} \cdot a^1).
\end{align*}
On the other hand, using 
lemma \ref{lambdaflemma} a) and the Takesaki-Takai isomorphism, we calculate 
\begin{align*}
\tr_A &(f \otimes (a^0 \rtimes r^0 \rtimes g^0) d(a^1 \rtimes r^1 \rtimes g^1)) \\
&= \tr_A(f(\hat{\psi}_{[2]} \sigma) \hat{\psi}_{[3]}(\delta) \F(\hat{\psi}_{[1]})
\otimes \gamma_A(a^0 \rtimes r^0 \rtimes g^0) d \gamma_A(a^1 \rtimes r^1 \rtimes g^1)) \\
&= f(\hat{\psi}_{[3]} \sigma) \hat{\psi}_{[4]}(\delta) \F(\hat{\psi}_{[2]})
\otimes S^{-1}(g^0)(\hat{\psi}_{(2)} S(r^1_{(4)})) S^{-1}(g^1)(\hat{\psi}_{(1)}) \\
&\qquad \psi(\hat{\psi}_{[1]} \hat{\psi}_{(3)} S(r^1_{(3)}) S(r^0_{(2)}) S^{-1}(\hat{\psi}_{(6)})) 
(\hat{\psi}_{(4)} S(r^1_{(2)}) \delta \hit S(r^0_{(1)}) \sigma \cdot a^0) \\
&\qquad d(\hat{\psi}_{(5)} \delta \hit S(r^1_{(1)}) \sigma \cdot a^1) \\
&= S^{-1}(g^0)(\hat{\psi}_{(2)} S(r^1_{(4)})) S^{-1}(g^1)(\hat{\psi}_{(1)})  
f(\hat{\psi}_{[3]} \sigma) \delta(\hat{\psi}_{[4]}) \F(\hat{\psi}_{[2]}) \\
&\qquad \psi(\hat{\psi}_{[1]} \hat{\psi}_{(3)} S(r^0_{(4)} r^1_{(3)}) S^{-1}(\hat{\psi}_{(6)})) \\
&\qquad \otimes (\hat{\psi}_{(4)} S(r^0_{(3)} r^1_{(2)}) \sigma \cdot a^0) 
d(\hat{\psi}_{(5)} \delta \hit S(r^0_{(2)} r^1_{(1)}) \sigma r^0_{(1)} \cdot a^1).
\end{align*}
Using $ \hat{\psi}(\delta h) = \hat{\psi}(h \delta) $ for all $ h \in \hat{H} $ this gives 
\begin{align*}
&\quad S^{-1}(g^0)(\hat{\psi}_{(2)} S(r^1_{(5)})) S^{-1}(g^1)(\hat{\psi}_{(1)}) f(\hat{\psi}_{[3]} \sigma) 
\F(\hat{\psi}_{[2]} \delta(r^0_{(4)} r^1_{(3)})) \\
&\qquad \psi(\hat{\psi}_{[1]} \hitby \delta \hat{\psi}_{(3)} \delta \hit S(r^0_{(5)} r^1_{(4)}) S^{-1}(\hat{\psi}_{(6)})) \\
&\qquad \otimes (\hat{\psi}_{(4)} S(r^0_{(3)} r^1_{(2)}) \sigma \cdot a^0)
d(\hat{\psi}_{(5)} \delta \hit S(r^0_{(2)} r^1_{(1)}) \sigma r^0_{(1)} \cdot a^1) \\
&= S^{-1}(g^0)(\hat{\psi}_{(2)} S(r^1_{(5)})) S^{-1}(g^1)(\hat{\psi}_{(1)}) f(\hat{\psi}_{[3]} \sigma) 
\F(\hat{\psi}_{[2]} \delta(r^0_{(4)} r^1_{(3)})) \\
&\qquad \psi(\hat{\psi}_{[1]} \hat{\psi}_{(3)} \hitby \delta^{-1} \sigma S^{-1}(r^0_{(5)} r^1_{(4)}) \sigma^{-1} S^{-1}(\hat{\psi}_{(6)})\hitby \delta^{-1}) \\
&\qquad \otimes (\hat{\psi}_{(4)} S(r^0_{(3)} r^1_{(2)}) \sigma \cdot a^0)
d(\hat{\psi}_{(5)} \delta \hit S(r^0_{(2)} r^1_{(1)}) \sigma r^0_{(1)} \cdot a^1)
\end{align*}
by invariance of $ \psi $ and $ \delta(\sigma) = 1 $. We get 
\begin{align*}
&\quad S^{-1}(g^0)(\hat{\psi}_{(2)} S(r^1_{(6)})) S^{-1}(g^1)(\hat{\psi}_{(1)}) f(\hat{\psi}_{[3]} \sigma) \\
&\qquad \F(\hat{\psi}_{[2]} \hat{\psi}_{(4)} \sigma S^{-1}(r^0_{(5)} r^1_{(4)}) \sigma^{-1} S^{-1}(\hat{\psi}_{(9)})
\hat{\psi}_{(8)} \sigma (r^0_{(4)} r^1_{(3)}) \hitby \delta \sigma^{-1} S(\hat{\psi}_{(5)})) \\
&\qquad \psi(\hat{\psi}_{[1]} \hat{\psi}_{(3)} \hitby \delta^{-1} \sigma S^{-1}(r^0_{(6)} r^1_{(5)}) \sigma^{-1} S^{-1}(\hat{\psi}_{(10)}) \hitby \delta^{-1}) \\
&\qquad \otimes (\hat{\psi}_{(6)} S(r^0_{(3)} r^1_{(2)}) \sigma \cdot a^0) 
d(\hat{\psi}_{(7)} \delta \hit S(r^0_{(2)} r^1_{(1)}) \sigma r^0_{(1)} \cdot a^1) \allowdisplaybreaks[2] \\
&= S^{-1}(g^0)(\hat{\psi}_{(2)} S(r^1_{(5)})) S^{-1}(g^1)(\hat{\psi}_{(1)}) f(\hat{\psi}_{[2]} \sigma) 
\F(\hat{\psi}_{(7)} \sigma (r^0_{(4)} r^1_{(3)}) \hitby \delta \sigma^{-1} S(\hat{\psi}_{(4)})) \\
&\qquad \psi(\hat{\psi}_{[1]} \hat{\psi}_{(3)} \hitby \delta^{-1} \sigma S^{-1}(r^0_{(5)} r^1_{(4)}) \sigma^{-1} S^{-1}(\hat{\psi}_{(8)}) \hitby \delta^{-1}) \\
&\qquad \otimes (\hat{\psi}_{(5)} S(r^0_{(3)} r^1_{(2)}) \sigma \cdot a^0)
d(\hat{\psi}_{(6)} \delta \hit S(r^0_{(2)} r^1_{(1)}) \sigma r^0_{(1)} \cdot a^1) 
\end{align*}
by right invariance of $ \psi $. This gives 
\begin{align*}
&\quad S^{-1}(g^0)(\hat{\psi}_{(2)} S(r^1_{(5)})) S^{-1}(g^1)(\hat{\psi}_{(1)}) f(\hat{\psi}_{[2]} \sigma) 
\F(\hat{\psi}_{(7)} \sigma (r^0_{(4)} r^1_{(3)}) \hitby \delta \sigma^{-1} S(\hat{\psi}_{(4)})) \\
&\qquad \psi(\hat{\psi}_{[1]} S^{-1}(\delta \hit \hat{\psi}_{(8)} \sigma r^0_{(5)} r^1_{(4)} \sigma^{-1} \delta \hit S(\hat{\psi}_{(3)})) \\
&\qquad \otimes (\hat{\psi}_{(5)} S(r^0_{(3)} r^1_{(2)}) \sigma \cdot a^0)
d(\hat{\psi}_{(6)} \delta \hit S(r^0_{(2)} r^1_{(1)}) \sigma r^0_{(1)} \cdot a^1) \allowdisplaybreaks[2] \\
&= S^{-1}(g^0)(\hat{\psi}_{(2)} S(r^1_{(4)})) S^{-1}(g^1)(\hat{\psi}_{(1)})  
f \bullet \F(\delta \hit \hat{\psi}_{(6)} \sigma (r^0_{(4)} r^1_{(3)}) \hitby \delta \sigma^{-1} \delta \hit S(\hat{\psi}_{(3)})) \\
&\qquad \otimes (\hat{\psi}_{(4)} S(r^0_{(3)} r^1_{(2)}) \sigma \cdot a^0) 
d(\hat{\psi}_{(5)} \delta \hit S(r^0_{(2)} r^1_{(1)}) \sigma r^0_{(1)} \cdot a^1)
\end{align*}
according to lemma \ref{lambdaflemma} b). Hence we get 
\begin{align*}
&\quad S^{-1}(g^0)(\hat{\psi}_{(2)} S(r^1_{(4)})) S^{-1}(g^1)(\hat{\psi}_{(1)})  
f \bullet \F(\delta \hit \hat{\psi}_{(6)} \delta \hit S^2(r^0_{(4)} r^1_{(3)}) \delta \hit S(\hat{\psi}_{(3)})) \\
&\qquad \otimes (\hat{\psi}_{(4)} S(r^0_{(3)} r^1_{(2)}) \sigma \cdot a^0) d(\hat{\psi}_{(5)} \delta \hit S(r^0_{(2)} r^1_{(1)}) \sigma r^0_{(1)} \cdot a^1) \\
&= S^{-1}(g^0)(\hat{\psi}_{(2)} S(r^1_{(5)})) S^{-1}(g^1)(\hat{\psi}_{(1)}) \delta(\hat{\psi}_{(7)}) \\
&\qquad f \bullet \F(\hat{\psi}_{(6)} \delta^{-1}(r^0_{(2)} r^1_{(1)}) \delta \hit S(\hat{\psi}_{(3)} S(r^0_{(5)} r^1_{(4)}))) \\
&\qquad \otimes (\hat{\psi}_{(4)} S(r^0_{(4)} r^1_{(3)})\sigma \cdot a^0) d(\hat{\psi}_{(5)} S(r^0_{(3)} r^1_{(2)}) \sigma r^0_{(1)} \cdot a^1)
\end{align*}
which shows that $ T^{-1} \tau_\ihd \tau_H $ and $ \tr_A $ agree on equivariant differential forms of degree one. \qed \\
In theorem \ref{technical} we have restricted ourselves to equivariant differential forms of degree zero and one. The assertion 
holds for higher degree forms as well, but the calculations become increasingly tedious to write down. For our 
purposes the above discussion is sufficient.

\section{Baaj-Skandalis duality and applications} \label{secbs}

In this section we use theorem \ref{technical} to formulate and prove the main result of this paper. As an application we derive 
the Green-Julg theorem in periodic cyclic homology for compact quantum groups, and the dual Green-Julg theorem in periodic 
cyclic cohomology for discrete quantum groups. \\ 
In the sequel we retain the general setup from the previous sections, in particular we assume that $ H $ is a bornological quantum 
group equipped with a modular pair. If $ B $ is an $ \hat{H} $-algebra then the map 
$ \tau_H: \Omega_H(\mathcal{T}(B) \rtimes \hat{H}) \rightarrow \Omega_\ihd(\mathcal{T}B) $ commutes with the boundary operators $ b_H $ and $ B_H $ 
according to proposition \ref{duality1}. Therefore it induces a chain map $ \tau_H: X_H(\mathcal{T}(B) \rtimes \hat{H}) \rightarrow X_\ihd(\mathcal{T}B) $. 
We define the duality morphism $ D_H: X_H(\mathcal{T}(B \rtimes \hat{H})) \rightarrow X_\ihd(\mathcal{T}B) $ by precomposing 
$ \tau_H $ with the chain map $ X_H(i_\ihd^B) $ induced by the canonical homomorphism 
$ i_\ihd^B: \mathcal{T}(B \rtimes \hat{H}) \rightarrow \mathcal{T}(B) \rtimes \hat{H} $. 
Now if $ A $ is an $ H $-algebra then the naturality of $ \tau $ implies
\begin{align*}
D_\ihd D_H = \tau_\ihd X_\ihd(i_H^A) &\tau_H X_H(i_\ihd^{A \rtimes H}) 
= \tau_\ihd \tau_H X_H(i_H^A \rtimes \hat{H}) X_H(i_\ihd^{A \rtimes H}).
\end{align*}
Taking into accout the Takesaki-Takai isomorphism $ A \rtimes H \rtimes \hat{H} \cong A \cotimes \mathcal{K}_H $, see 
proposition \ref{TTduality},  
we conclude that $ D_\ihd D_H: X_H(\mathcal{T}(A \rtimes H \rtimes \hat{H})) \rightarrow X_\ihd(\mathcal{T}A) $ 
coincides with the map $ t_A $ defined in the proof of theorem 8.4 in \cite{Voigtechq}. 
Therefore theorem \ref{technical} implies together with stability of equivariant periodic cyclic homology \cite{Voigtechq} that the chain map 
\begin{equation*}
D_\ihd D_H: X_H(\mathcal{T}(A \rtimes H \rtimes \hat{H} \rtimes H \rtimes \hat{H}))
\rightarrow X_H(\mathcal{T}(A \rtimes H \rtimes \hat{H}))
\end{equation*}
is a homotopy equivalence of pro-paracomplexes of Yetter-Drinfeld modules for every $ H $-algebra $ A $. \\
As a consequence, we obtain the following duality result.
\begin{theorem} \label{BaajSkandalisduality2}
Let $ H $ be a bornological quantum group equipped with a modular pair. Then there exists a natural isomorphism 
$$
J_H: HP^H_*(A,B) \rightarrow HP^{\hat{H}}_*(A \rtimes H, B \rtimes H)
$$
for all $ H $-algebras $ A $ and $ B $. This isomorphism is compatible with 
composition products and maps the class of an $ H $-equivariant homomorphism $ A \rightarrow B $ 
to the class of the induced homomorphism $ A \rtimes H \rightarrow B \rtimes H $. 
\end{theorem}
\proof We may identify $ HP^H_*(A,B) \cong HP^H_*(A \rtimes H \rtimes \hat{H}, B \rtimes H \rtimes \hat{H}) $ according to 
stability and Takesaki-Takai duality, and we define $ J_H $ by 
$$
J_H(\phi) = D_H \circ \phi \circ D_H^{-1}, 
$$
using that $ D_H: X_H(\mathcal{T}(A \rtimes H \rtimes \hat{H} \rtimes H \rtimes \hat{H})) \rightarrow 
X_\ihd(\mathcal{T}(A \rtimes H \rtimes \hat{H} \rtimes H)) $ is a homotopy equivalence, see above. 
By naturality of $ D_H $ it follows that $ J_H $ is natural and satisfies $ J_H([f]) = [f \rtimes H] $ if $ f: A \rightarrow B $ 
is an $ H $-equivariant homomorphism and $ f \rtimes H: A \rtimes H \rightarrow B \rtimes H $ the corresponding homomorphism 
of the crossed products. Moreover it is immediate from the definitions that 
$ J_H $ is multiplicative. \qed \\
Let us now derive a general version of the 
Green-Julg theorem in cyclic homology for compact quantum groups, compare \cite{Brylinskibrown}, \cite{Blockthesis}, \cite{Buesthesis}, \cite{AK2}. 
If $ G $ is a compact quantum group we write $ \Pol(G) $ for the unital Hopf $ * $-algebra of polynomial functions on $ G $. 
\begin{theorem} \label{greenjulg}
Let $ H $ be the bornological quantum group dual to the polynomial algebra $ \Pol(G) $ of a compact quantum group $ G $. Then there is a natural isomorphism 
$$
HP^H_*(\mathbb{C},A) \cong HP_*(A \rtimes H) 
$$
for all $ H $-algebras $ A $. 
\end{theorem}
\proof We consider the canonical modular pair $ (1, \delta) $ for $ \Pol(G) $, where $ \delta = f_{1/2} $ is the modular character. 
The left hand side of the asserted isomorphism is computed by the homology of 
$ X_H(\T(A \rtimes H \rtimes \hat{H}))^H $, that is, the $ H $-invariant part of the equivariant $ X $-complex $ X_H(\T(A \rtimes H \rtimes \hat{H})) $, 
taken on each level of the underlying inverse system.
From our above computations and stability it follows that this complex is 
homotopy equivalent to $ X_\ihd(\T(A \rtimes H))^H $. Remark that the Yetter-Drinfeld module action of $ H $  
coincides with the anti-Yetter-Drinfeld module action in this situation since the group-like element of our modular pair 
is equal to $ 1 $. Moreover, is is easy to check that an element of $ \Omega_\ihd(A \rtimes H) $ is invariant under the action of $ H $ iff it 
is of the form $ \hat{\phi} \otimes \omega $ for some $ \omega \in \Omega(A \rtimes H) $ where $ \hat{\phi}: \Pol(G) \rightarrow \mathbb{C} $ is 
the left and right invariant Haar integral. This yields 
an isomorphism $ X_\ihd(\T(A \rtimes H))^H \cong X(\T(A \rtimes H)) $. 
It is easy to check that this isomorphism is compatible with the boundary operators, thus finishing the proof. \qed \\
Another application of our computations is the following dual Green-Julg theorem for discrete quantum groups. 
We recall that a discrete quantum group $ G $ can be viewed as the dual of a compact quantum group, 
and we write $ \mathbb{C}[G] $ for the associated unital Hopf $ * $-algebra in this case. 
\begin{theorem} \label{dualgreenjulg}
Let $ H = \mathbb{C}[G] $ be the polynomial group algebra of a discrete quantum group $ G $. Then there is a natural isomorphism 
$$
HP^H_*(A, \mathbb{C}) \cong HP^*(A \rtimes H) 
$$
for all $ H $-algebras $ A $. 
\end{theorem}
\proof Again we let $ (1, \delta) $ be the canonical modular pair for $ H $. In the same way 
as in the proof of theorem 16.4 in \cite{Voigtepch} we see that the left hand side of the 
asserted isomorphism is computed by $ H_*(\Hom_H(X_H(\T(A \rtimes H \rtimes \hat{H})), \mathbb{C}_\delta[0])) $, where 
$ \mathbb{C}_\delta[0] $ denotes the complex numbers viewed as a trivial complex in degree zero 
with the $ H $-action given by $ t \cdot 1 = \delta^{-1}(t) 1 $. \\
We thus have to study the dual space of $ X_H(\T(A \rtimes H \rtimes \hat{H}))_\delta $, the latter denoting the twisted coinvariants 
of $ X_H(\T(A \rtimes H \rtimes \hat{H})) $. The twisted coinvariants are obtained by taking 
the quotient by the closed linear span of all elements of the form 
$$ 
t \bullet m - \delta^{-1}(t) m 
$$ 
for $ t \in H $ in each degree of the inverse system of Yetter-Drinfeld modules underlying $ X_H(\T(A \rtimes H \rtimes \hat{H})) $. 
Using our above computations, the same arguments as in section 16 of \cite{Voigtepch} show that the 
duality map $ D_H $ induces a homotopy equivalence between $ X_H(\T(A \rtimes H \rtimes \hat{H}))_\delta $ and $ X_\ihd(\T(A \rtimes H))_\delta $. \\ 
Let us define a bounded linear map 
$ \alpha: X_\ihd(\T(A \rtimes H))_\delta \rightarrow X(\T(A \rtimes H)) $ by 
$ \alpha(t \otimes \omega) = \delta^{-1}(t) \omega $. First observe that $ \alpha $ is well-defined 
since 
$$
\alpha(t \bullet (r \otimes \omega)) = \alpha(tr \otimes \omega) = \delta^{-1}(tr) \omega = \alpha(\delta^{-1}(t) r \otimes \omega). 
$$
Similarly, we define $ \beta: X(\T(A \rtimes H)) \rightarrow 
X_\ihd(\T(A \rtimes H))_\delta $ by $ \beta(\omega) = 1 \otimes \omega $ where $ 1 \in H $ is the unit element. 
It is easy to see that $ \alpha $ and $ \beta $ are inverse isomorphisms compatible with the boundary operators. 
This proves the claim. \qed \\
We remark that theorem \ref{greenjulg} and theorem \ref{dualgreenjulg} can also be derived more formally from theorem 
\ref{BaajSkandalisduality2}, by basically reducing everything to statements abouts adjoint functors. 
However, the arguments above make it easier to keep track of the maps on the chain level.

\bibliographystyle{plain}

\bibliography{cvoigt}

\end{document}